\documentclass[a4paper,11pt]{article}
\usepackage{titlesec}
\usepackage{amssymb,amsthm,amsmath}   % AMS 数学公式的
\usepackage{cases}
\usepackage{cite}
\usepackage{bm}
\usepackage{color}    % 可以随意定义文中字体颜色
\usepackage[dvips]{graphicx}
\usepackage{subfigure}
\usepackage{tikz}
\usetikzlibrary{shapes.geometric} % 可以畫菱形等非矩形、圓的 shape
\usetikzlibrary{calc} % 座標計算

\newtheorem{example}{Example}

\tikzset{
passprocess/.style={rectangle, minimum width=75pt, minimum height=25pt, draw=black},%, font=\ttfamily, text centered},
startstop/.style={rectangle, rounded corners=5pt, minimum width=75pt, minimum height=25pt, draw=black},%, font=\ttfamily, text centered, fill=orange},
decision/.style={
diamond,% 菱形
shape aspect=3,%aspect value is the ratio of width and height for diamond
minimum width=75,
draw=black},%,%the color of line
line/.style={draw, ->, %shorten>=2pt,
thick
}
} % end of tikzset.

\begin{document}
% =====================================================
% =                标题与作者                           =
% =====================================================
\title{\bfseries Numerical investigations of traveling singular sources problems via moving mesh method % in conjunction with domain decomposition
  % A novel moving mesh method based on the domain decomposition for
  % traveling singular sources problems
}
\renewcommand{\thefootnote}{}
\author{Zhicheng Hu$^{a,}$\thanks{Corresponding author.% \emph{E-mail address:} huzhicheng1986@gmail.com (Z. Hu); matlkw@zju.edu.cn (K. Liang)
}\footnotemark[1], ~Keiwei Liang$^b$
    \\ \small \emph{$^a$LMAM \emph{\&} School of Mathematical Sciences, Peking University, Beijing 100871, China}
    \\ \small \emph{$^b$Department of Mathematics, Zhejiang University, Hangzhou 310027, China}
}

\date{}   % 取消显示日期

\maketitle   % 打出标题
% 貌似这样使用脚注的时候要放到 \maketitle 后面, 否则这个脚注打不出来...
\footnotetext[1]{\emph{E-mail address:} huzhicheng1986@gmail.com (Z. Hu); matlkw@zju.edu.cn (K. Liang).}

%\vspace{-1em}

% ===============================================
% =                摘要                          =
% ===============================================
\begin{abstract}
  This paper studies the numerical solution of traveling singular
  sources problems. In such problems, a big challenge is the sources
  move with different speeds, which are described by some ordinary
  differential equations. A predictor-corrector algorithm is presented
  to simulate the position of singular sources. Then a moving mesh
  method in conjunction with domain decomposition is derived for the
  underlying PDE. According to the positions of the sources, the whole
  domain is splitted into several subdomains, where moving mesh
  equations are solved respectively. On the resulting mesh, the
  computation of jump $[\dot{u}]$ is avoided and the discretization of
  the underlying PDE is reduced into only two cases. In addition, the
  new method has a desired second-order of the spatial
  convergence. Numerical examples are presented to illustrate the
  convergence rates and the efficiency of the method. Blow-up
  phenomenon is also investigated for various motions of the sources.

  \vspace{1em} \noindent{\bf Keywords:} Moving mesh method; Domain
  decomposition; Traveling singular sources
\end{abstract}

% \vspace{1em}

% ==============================================
% =                 正文                        =
% ==============================================
% 第一节
\section{Introduction}
We take the one-dimensional moving singular sources equation
\begin{align}
  \label{eq:heat_eq}
  & u_t - u_{xx} = \sum_{i=0}^{q-1}F_i(t,x,u) \delta(x-\alpha_i(t)), \quad -\infty < x < \infty,~ t>0, \\
  \label{eq:init_cond}
  & u(x,0) = u_0(x), \quad -\infty <x<\infty,\\
  \label{eq:bound_cond}
  & u(x,t) \to 0 \quad \text{as} \quad |x| \to \infty, ~t>0.
\end{align}
as the model problem in this paper. Here $q>0$ is the number of singular
sources.  The initial value $u_0(x)$ is taken to be continuous and
compatible with the boundary conditions, i.e. $u_0(x)\to 0$ as $|x|\to
\infty$. The local source functions $F_i(t, x, u)$
$(i=0,1,\ldots,q-1)$ might be given a priori or can be determined from
some additional constraints on the solution.  The traveling sources
are located at $\alpha_i(t)$, $i=0,1,\ldots,q-1$.  In general, their
velocities can be described by several ordinary differential equations
\begin{align}
  \label{eq:source_eq}
  \frac{\mathrm{d} \alpha_i}{\mathrm{d} t} = \psi_i ( t, \alpha_i(t), u ), \quad i = 0, 1, \ldots, q-1,
\end{align}
which are coupled with the solution $u$. We assume that the sources do
not intersect with each other during the time in consideration. This
model arises in many areas such as laser beams traveling problems
where $u$ is the temperature of the material \cite{kirk2002blow}, or
free-boundary solidification problems where $\alpha_i(t)$,
$i=0,1,\ldots,q-1$, are the moving interfaces between different phases
\cite{beyer1992analysis}.

It is well known that the solution of the model is continuous and
piecewise smooth \cite{li1997immersed}.  However, the derivative of
the solution has a jump at each source due to the delta function
singularity on it, and the jump is given by \cite{ma2009moving}
\begin{align}
  \label{eq:jump_ux}
  [u_x]_{(\alpha_i(t), t)} = -F_i(t, \alpha_i(t), u(\alpha_i(t),t)), \quad \quad i = 0, 1, \ldots, q-1.
\end{align}
This leads the standard numerical methods, either finite difference
method or finite element method, might fail when crossing the
time-dependent source positions. Various approaches, have been used to
deal with the delta function singularity, such as the immerse boundary
(IB) method and the immerse interface method (IIM)
\cite{smereka2006numerical,beyer1992analysis,leveque1994immersed,
  li1997immersed, kandilarov2003immersed, yang2009smoothing}. For the
IB method originally proposed by \cite{peskin1977numerical}, the delta
function is approximated by an appropriately chosen discrete delta
function.  Beyer and LeVeque \cite{beyer1992analysis} studied various
cases of the model \eqref{eq:heat_eq}$-$\eqref{eq:bound_cond} for the
IB method with $q=1$, and the source position $\alpha_i(t)$ being
priori specified.  In contrast, the IIM first introduced by LeVeque
and Li \cite{leveque1994immersed} incorporates the known jumps of
solution or its derivatives into the finite difference scheme to
obtain a modified discretization scheme. Li \cite{li1997immersed}
developed an IIM numerical algorithm on the uniform mesh for the model
\eqref{eq:heat_eq}$-$\eqref{eq:bound_cond} and \eqref{eq:source_eq}
with $q=1$.

The model \eqref{eq:heat_eq}$-$\eqref{eq:bound_cond} and
\eqref{eq:source_eq} is more difficult to be solved when the source
function $F_i(t,x,u)$ is high nonlinearity.  In this case, the
solution is always blow-up in some finite time $T>0$ if the sources
are stationary or move at sufficiently low speed, while blow-up will
be avoided if the sources move at sufficiently high speed (see
e.g. \cite{olmstead1997critical,kirk2000influence,kirk2002blow}).  As
the solution evolves singularity, the uniform mesh method always
become computationally prohibitive. Hence, moving mesh method has to
be employed, %we should employ an adaptive mesh method.
which is one of the most popular adaptive methods and have been
successfully used to investigate the blow-up phenomenon
\cite{budd1996moving,huang2008study}.  In MMPDE's approaches, the
movement of the mesh is controlled by the moving mesh partial
differential equations (MMPDEs) based on the equidistribution
principle \cite{huang1994MMPDE}.  Among these MMPDEs, MMPDE4, MMPDE5
and MMPDE6 are popular to use.  Readers interested in the moving mesh
method and its applications can refer to the books
\cite{tang2007adaptive,Huang2011book}.

Recently, several papers have been devoted to moving mesh method for
the model \eqref{eq:heat_eq}$-$\eqref{eq:bound_cond} with a priori
specified source position $\alpha_i(t)$ for $q=1$
\cite{zhu2010numerical,ma2009moving} and for $q>1$, in which the
sources move with the same speed \cite{zhu_moving,hu2011moving}.  This
paper is a further study of \cite{hu2011moving} and \cite{hu2012zh}
for the model \eqref{eq:heat_eq}$-$\eqref{eq:bound_cond} with general
movement of the sources, which do not intersect with each other during
time evolution. First, we choose a finite observed domain containing
all sources with appropriate boundary conditions, and divide it into
$q+1$ subdomains by the $q$ sources. Obviously, the sizes of the
subdomains are changed as the sources traveling. MMPDEs are applied on
each subdomain to obtain a local equidistributed mesh on it. Then the
underlying PDE \eqref{eq:heat_eq} is solved on the whole observed
domain with the mesh composed of the local mesh on each
subdomain. Taking the advantages of domain decomposition
\cite{Toselli2005book}, MMPDEs could be solved efficiently by parallel
computing. Moreover, It can be found that the computation of
$[\dot{u}]$ is avoided, thus the discretization scheme for the
underlying PDE becomes very simple. In addition, our method has an
expected second-order convergence in space.

The organization of the paper is as follows. In section
\ref{sec:DDMMM}, we introduce the moving mesh method in conjunction
with domain decomposition for the model problem. In section
\ref{sec:discretization}, the discretization schemes for the physical
problem will be derived in detail. In section \ref{sec:example},
several numerical examples are given to demonstrate the numerical
efficiency and accuracy of our method. The conclusions are presented
in the last section.%\ref{sec:conclusion}.

\section{Moving mesh method in conjunction with domain decomposition}\label{sec:DDMMM}   % 其他地方要引用的时候用 \ref{sec1}
In the last decade, moving mesh method in conjunction with a Schwarz
domain decomposition has been developed by Haynes and his co-workers
(see i.e. \cite{haynes2007schwarz, gander2012domain} and references
therein). And in this section, we will introduce a slight different
method, that is, moving mesh method in conjunction with a
non-overlapping domain decomposition. % we first introduce the domain
% decomposition by the sources, and then briefly describe the moving
% mesh strategy used in our method.

%    \subsection{Domain decomposed}
Denote the observed domain by $[x_l,x_r]$ and assume it containing all
sources, that is,
$x_l<\alpha_0(t)<\alpha_1(t)<\cdots<\alpha_{q-1}(t)<x_r$. Here $x_l$,
$x_r$ are either constants or variables of $t$. Then, the observed
domain is divided into $q+1$ subdomains $[\alpha_{i-1},\alpha_{i}]$
$(i=0,1,\ldots,q)$ with $\alpha_{-1}=x_l$, $\alpha_{q}=x_r$, by the
$q$ sources respectively. Obviously, the sizes of the subdomains are
variables of $t$ too.

Let $x$ and $\xi$ denote physical and computational coordinates,
respectively. Without loss of generality we assume the computational
domain is $[0,1]$. Then an one-to-one coordinate transformation
between the observed domain $[x_l,x_r]$ and the computational domain
$[0,1]$ is defined by
    \begin{equation}\label{eq:logical_map}
    x = x(\xi,t), \quad \xi \in [0,1],
    \end{equation}
with
    \begin{equation*}
    x(0,t)=x_l, \quad x(1,t) = x_r.
    \end{equation*}
    For a given uniform mesh, $\xi_j=\frac{j}{N}$, $j=0,1,\ldots,N$,
    on the computational domain, the corresponding mesh on the
    observed domain $[x_l,x_r]$ is
    \begin{equation*}
    x_l=x_0(t)<x_1(t)<\cdots<x_{N-1}(t)<x_N(t)=x_r.
    \end{equation*}
    In our method, the coordinate transformation
    \eqref{eq:logical_map} is determined as a piecewise smooth
    function. On each subdomain $[\alpha_{i-1},\alpha_{i}]$,
    $i=0,1,\ldots,q$, it is the solution of an MMPDE which is derived
    from the equidistribution principle. In the literature, the
    following MMPDEs
    % Following \cite{huang1994MMPDE}, MMPDE4, MMPDE5, and MMPDE6 are
    % given in the following
    \begin{equation}\label{eq:MMPDE4}
    \frac{\partial}{\partial\xi} \bigg(M \frac{\partial{\dot{x}}}{\partial\xi}\bigg)
    =-\frac{1}{\tau}\frac{\partial}{\partial\xi}
    \bigg(M\frac{\partial x}{\partial \xi}\bigg),
    \end{equation}
    \begin{equation}\label{eq:MMPDE5}
    -\dot{x}=-\frac{1}{\tau}\frac{\partial}{\partial\xi}
    \bigg(M\frac{\partial x}{\partial \xi}\bigg),
    \end{equation}
    \begin{equation}\label{eq:MMPDE6}
    \frac{\partial^2\dot{x}}{\partial\xi^2}=-\frac{1}{\tau}\frac{\partial}{\partial\xi}
    \bigg(M\frac{\partial x}{\partial \xi}\bigg),
    \end{equation}
    which known as MMPDE4, MMPDE5 and MMPDE6, respectively, are
    popularly used after they were originally established and analyzed
    in \cite{huang1994MMPDE}. Here $M=M(x,t)$ is the monitor function
    giving some measure of the solution error on the physical domain
    and $\tau>0$ is a parameter representing a timescale for adjusting
    the mesh toward equidistribution. In the asymptotic case
    $t\to\infty$, the solution of MMPDE4, MMPDE5 and MMPDE6 would
    satisfy the equidistribution principle, which is stated that
    \cite{huang1994MMPDE}
    \begin{align}
      \label{eq:EP}
\frac{\partial}{\partial\xi}
\bigg(M\frac{\partial x}{\partial \xi}\bigg) = 0.
    \end{align}
    For more details about MMPDE, one can refer to
    \cite{huang1994MMPDE} or the recent book \cite{Huang2011book}. In
    this paper, MMPDE6 \eqref{eq:MMPDE6} with the boundary condition
\begin{align}
  \label{eq:BV_MMPDE}
  x(\xi_{j_{i-1}^s},t) = \alpha_{i-1}(t), \quad x(\xi_{j_i^s},t) = \alpha_i(t), \qquad i=0,1,\ldots,q,
\end{align}
is employed as an example to describe our moving mesh strategy in
conjunction with domain decomposition. Here $j_i^s$ is some fixed
index satisfying $0<j^s_i<N$. The resulting mesh, used to solve the
model problem on $[x_l,x_r]$, satisfies the property that a fixed mesh
point is located on each source during the time in consideration,
i.e. $x_{j^s_i}\equiv \alpha_i(t)$.

\begin{figure}[!htbp]
  \centering
    \small
%\begin{minipage}{0.75\textwidth}
  \begin{tikzpicture}[scale = 0.75]
    \matrix[row sep=5mm,column sep=5mm, top color=white, bottom color=white, scale=0.75]{
      \node[startstop, text width=0.5\textwidth, scale=0.75, text centered] (start) {Given the old mesh $x_j^n$ on the observed domain $[x_l,x_r]$ and the corresponding solution on the mesh.};\\
      \node[passprocess, scale=0.75] (monitor) {Compute the monitor function $M$ on the mesh.};\\
      \node[passprocess, text width=0.5\textwidth, scale=0.75] (MMPDE) {Solve MMPDE6 \eqref{eq:MMPDE6} with the boundary condition \eqref{eq:BV_MMPDE} on each subdomain $[\alpha_{i-1},\alpha_i]$, $i=0,1,\ldots,q$.};\\
      \node[startstop, text width=0.5\textwidth, scale=0.75, text centered] (end) {Combining the local equidistributed mesh on each subdomain to give the new mesh $x_j^{n+1}$ on $[x_l,x_r]$.}; \\
    };
  \begin{scope}[every path/.style=line]
    \path (start) -- (monitor);
    \path (monitor) -- (MMPDE);
    \path (MMPDE) -- (end);
  \end{scope}
  \end{tikzpicture}
%\end{minipage}
  \caption{The moving mesh strategy in conjunction with domain
    decomposition.}
  \label{fig:mesh_strategy}
\end{figure}
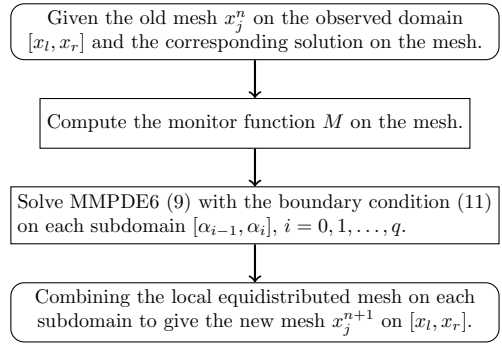

Figure \ref{fig:mesh_strategy} shows the moving mesh strategy in
conjunction with domain decomposition. Here the computation of the
monitor function will be presented in section \ref{sec:example}. And
MMPDE6 \eqref{eq:MMPDE6} is solved by the following finite difference
scheme
    \begin{multline}\label{eq:MMPDE6_dis}
    \frac{\big(x_{j+1}^{n+1}-2x_j^{n+1}+x_{j-1}^{n+1}\big) -
    \big(x_{j+1}^n-2x_j^n+x_{j-1}^n\big)}{\Delta t_n} \\
    =
    -\frac{1}{\tau}\bigg(M_{j+\frac{1}{2}}\big(x_{j+1}^{n+1}-x_j^{n+1}\big)
    -
    M_{j-\frac{1}{2}}\big(x_{j}^{n+1}-x_{j-1}^{n+1}\big)\bigg)
    \end{multline}
    in our numerical examples, where $\Delta t_n = t_{n+1}-t_{n}$ and
    $M_{j+\frac{1}{2}}=(M_{j+1}+M_{j})/2$.

% We conclude this section with some remarks. On the one hand, t
    The new mesh could be obtained very efficiently by parallel
    computing based on domain decomposition methods
    \cite{Toselli2005book}. And it is best in the sense of
    equidistribution on each subdomain.  On the other hand, we will
    found in the next section that the computation of the jump
    $[\dot{u}]$ is avoided, hence the discretization scheme for
    the physical PDE \eqref{eq:heat_eq} becomes very simple.

\section{Model discretization and final algorithm}\label{sec:discretization}
In this section, we derive the discretization schemes for the physical
model problem \eqref{eq:heat_eq}$-$\eqref{eq:bound_cond} and
\eqref{eq:source_eq} on the observed domain $[x_l,x_r]$ with
appropriate boundary conditions. Then present a full algorithm of
moving mesh method for the model problem.

\subsection{Discretization schemes}
For an arbitrary function $f=f(x,t)=f(x(\xi,t),t)$, we have
    \begin{equation*}
        \dot{f}=\frac{\partial{f}}{\partial{t}}(x(\xi,t),t)
            \bigg|_{\xi~\text{fixed}}
         = f_{t} + f_{x}\dot{x}.
    \end{equation*}
    Through the coordinate transformation \eqref{eq:logical_map}, we
    can rewrite equation \eqref{eq:heat_eq} on the computational
    coordinates as
\begin{align}
  \label{eq:heat_logical}
  \dot{u} - u_x \dot{x} - u_{xx} = \sum_{i=0}^{q-1}F_i(t,x,u) \delta(x-\alpha_i(t)).
\end{align}
Since the right-hand side of \eqref{eq:heat_logical} vanishes when
$x\neq \alpha_i(t)$, that is,
    \begin{equation}\label{eq:reduce_logical}
      \dot{u} - u_x \dot{x} - u_{xx} = 0,
    \end{equation}
    we conduct the discretization scheme for \eqref{eq:heat_logical}
    on the above equation as \cite{hu2011moving}, with each term on
    the left-hand side of \eqref{eq:reduce_logical} containing the
    information of jumps when they cross the sources. Physically, the
    value $u$ of $i$th source changes smoothly as time evolution,
    which means the jump of the directional derivative of $u(x,t)$
    along the vector $(\alpha_i'(t),1)$ is zero \cite{ma2009moving,
      hu2011moving}, i.e.,
    % Following \cite{ma2009moving} and \cite{hu2011moving}, we have
    % that the jump of the directional derivative of $u(x,t)$ along
    % the vector $(\alpha_i'(t),1)$ is zero, that is,
    \begin{equation}\label{eq:jump_u_t}
      [u_t]_{(\alpha_i(t),t)}+[u_x]_{(\alpha_i(t),t)}\alpha_i'(t)=0, \quad i=0,1,\ldots,q-1.
    \end{equation}
Recalling that $x_{j^s_i}\equiv \alpha_i(t)$, it follows that
    \begin{equation}\label{eq:jump_u_dot}
    [\dot{u}]_{(\alpha_i(t),t)} =
    [u_{t}+u_{x}\dot{x}]_{(\alpha_i(t),t)}=
    [u_{t}]_{(\alpha_i(t),t)}+\dot{x}[u_{x}]_{(\alpha_i(t),t)}=0,
\quad i=0,1,\ldots,q-1.
    \end{equation}
    By using the above equation, we can deduce from
    \eqref{eq:reduce_logical} that
    \begin{equation}\label{eq:jump_u_xx_tmp}
      [u_{xx}]_{(\alpha_i(t),t)} = [\dot{u}-u_x\dot{x}]_{(\alpha_i(t),t)} =
      [\dot{u}]_{(\alpha_i(t),t)} - \dot{x} [u_{x}]_{(\alpha_i(t),t)} =
      - \alpha_i'(t) [u_{x}]_{(\alpha_i(t),t)},
      % = \psi_i(t,\alpha_i(t),u)F_i(t,\alpha_i(t),u(\alpha_i(t),t)),
    \end{equation}
    $i=0,1,\ldots,q-1$. Then we obtain immediately
    \begin{equation}\label{eq:jump_u_xx}
      [u_{xx}]_{(\alpha_i(t),t)} = \psi_i(t,\alpha_i(t),u)F_i(t,\alpha_i(t),u(\alpha_i(t),t)), \quad i=0,1,\ldots,q-1,
    \end{equation}
    by taking \eqref{eq:source_eq} and \eqref{eq:jump_ux} into
    \eqref{eq:jump_u_xx_tmp}.
%     Combining with \eqref{eq:heat_eq}, \eqref{eq:source_eq},
%     \eqref{eq:jump_ux} and \eqref{eq:jump_u_t}, we have that
%     \begin{equation}\label{eq:jump_u_xx}
%     [u_{xx}]_{(\alpha_i(t),t)} =
%     [u_{t}]_{(\alpha_i(t),t)}=
%     -\alpha_i'(t)[u_{x}]_{(\alpha_i(t),t)}=\psi_i(t,\alpha_i(t),u)F_i(t,\alpha_i(t),u(\alpha_i(t),t)),
%     \end{equation}
% $i=0,1,\ldots,q-1$.

Similarly to \cite{hu2011moving}, the discretization scheme for
\eqref{eq:heat_logical} are divided into two cases due to
$x_{j^s_i}\equiv \alpha_i(t)$ during time integration. For $j\neq
j^s_i$, $i=0,1,\ldots,q-1$, the mesh point not located at the source,
\eqref{eq:heat_logical} is discretized by standard center difference
for spatial variable and backward difference for temporal variable,
that is,
    \begin{equation}\label{eq:phy_dis_regular}
    \frac{u_j^{n+1}-u_j^{n}}{\Delta t_n} -
    \frac{u_{j+1}^{n+1}-u_{j-1}^{n+1}}{h_{j+1}^{n+1}+h_{j}^{n+1}}
    \bigg(\frac{x_j^{n+1}-x_j^n}{\Delta t_n}\bigg) - \frac{2}{h_{j+1}^{n+1}+h_{j}^{n+1}} \bigg(
    \frac{u_{j+1}^{n+1}-u_j^{n+1}}{h_{j+1}^{n+1}} -
    \frac{u_{j}^{n+1}-u_{j-1}^{n+1}}{h_{j}^{n+1}} \bigg) = 0,
    \end{equation}
    where $h_j^n = x_j^{n}-x_{j-1}^n$. Here $x_j^n$, $u_j^n$ are the
    mesh and the solution on it at time step $t_n$, respectively. For
    $j=j^s_i$, the mesh point just located at the source, the jump
    informations should be incorporated into the discretization
    scheme. For this case, the discretization scheme for
    \eqref{eq:heat_logical} reads
\begin{align}\label{eq:phy_dis_irregular}
  \frac{u_{j^s_i}^{n+1}-u_{j^s_i}^n}{\Delta t_n} -
  \frac{u_{j^s_i+1}^{n+1}-u_{j^s_i-1}^{n+1}}
  {h_{j^s_i+1}^{n+1}+h_{j^s_i}^{n+1}} \psi_i^{n+1} & -
  \frac{2}{h_{j^s_i+1}^{n+1}+h_{j^s_i}^{n+1}} \bigg(
  \frac{u_{j^s_i+1}^{n+1}-u_{j^s_i}^{n+1}}{h_{j^s_i+1}^{n+1}} -
  \frac{u_{j^s_i}^{n+1}-u_{j^s_i-1}^{n+1}}{h_{j^s_i}^{n+1}} \bigg)\nonumber\\
  & - \frac{2}{h_{j^s_i+1}^{n+1}+h_{j^s_i}^{n+1}} F_i(u_{j^s_i}^{n+1})
  = 0,
\end{align}
where $\psi_i^{n+1}\approx \psi_i(t_{n+1}, \alpha_i^{n+1}, u^{n+1}) $,
$F_i(u_{j^s_i}^{n+1})\approx F_i(t_{n+1}, \alpha_i^{n+1},
u_{j^s_i}^{n+1})$, $i=0,1,\ldots,q-1$.

In the above schemes, we need the source position $\alpha_i(t_{n+1})$
at time step $t_{n+1}$. For the general movement \eqref{eq:source_eq},
it is computed by the following Crank-Nicolson scheme
\begin{align}
  \label{eq:source_CN}
  \alpha_i^{n+1} = \alpha_i^n + \frac{\Delta t_n}{2} (\psi_i^{n+1}+\psi_i^{n}), \quad i=0,1,\ldots,q-1,
\end{align}
as in \cite{hu2012zh}. If $\psi_i(t,\alpha_i(t),u)$,
$i=0,1,\ldots,q-1$, are independent of $u$, the source position
$\alpha_i^{n+1}$ and the speed $\psi_i^{n+1}$ can be calculated in
advance before solving the discretization schemes for MMPDE6
\eqref{eq:MMPDE6} and physical PDE \eqref{eq:heat_logical}. Otherwise,
the resulting system would be too complicated to be solved. In this
case, we decouple the discretization system by a predictor-corrector
algorithm. For the predictor step, assume $\psi_i^{n+1}=\psi_i^n$ and
solve \eqref{eq:source_CN} to get an approximate variable $\alpha_i^*$
of $\alpha_i^{n+1}$. Then substituting $\psi_i^{n+1}$ and $\alpha_i^*$
into the discretization schemes for \eqref{eq:MMPDE6} and
\eqref{eq:heat_logical} to obtain an approximate solution $u^*$ of
$u^{n+1}$. For the corrector step, compute
$\psi_i^{n+1}=\psi_i(t_{n+1},\alpha_i^*,u^*)$ and solve
\eqref{eq:source_CN} to get $\alpha_i^{n+1}$, then obtain the solution
$u^{n+1}$ at time step $t_{n+1}$ by the discretization schemes for
\eqref{eq:MMPDE6} and \eqref{eq:heat_logical}.

To complete the discretization schemes, we require an appropriate
condition for $u$ on the boundary of the observed domain $[x_l,x_r]$.
For the observed domain is small enough, we employ a third-order local
absorbing boundary condition (LABC) proposed in
\cite{brunner2010computational}
\begin{equation}\label{eq:LABCs_phy}
  3s_0 u_x + u_{xt} \pm s_0\sqrt{s_0}u\pm
  3\sqrt{s_0}u_t = 0
\end{equation}
for \eqref{eq:heat_eq} as in \cite{hu2011moving,hu2012zh}.  Here $s_0$
is an user-defined parameter, the plus sign in "$\pm$" corresponds to
the LABC at the right boundary $x_r$, and the minus sign corresponds
to the one at the left boundary $x_l$. Under the map
\eqref{eq:logical_map}, we get the LABC for \eqref{eq:heat_logical} as
follows
\begin{equation}\label{eq:LABCs_log}
  \dot{u}_{x} \pm 3 \sqrt{s_0} \dot{u} + 3 s_0 u_{x}-
  u_{xx}\dot{x} \pm s_0\sqrt{s_0}u \pm \left( - 3\sqrt{s_0}u_{x}\dot{x}\right) = 0,
\end{equation}
where the plus sign in "$\pm$" corresponds to the right boundary, and
the minus sign corresponds to the left boundary. According to
\cite{hu2011moving,hu2012zh}, a finite difference scheme for
\eqref{eq:LABCs_log} is
\begin{align}\label{eq:LABC_left}
  & \frac{1}{\Delta t_n}\bigg(\frac{u_1^{n+1}-u_{-1}^{n+1}}{2h_1^{n+1}} -
  \frac{u_1^{n}-u_{-1}^{n}}{2h_1^{n}}\bigg) - 3 \sqrt{s_0}
  \frac{u_0^{n+1}-u_0^{n}}{\Delta t_n} + 3\left(s_0+\sqrt{s_0}\frac{x_0^{n+1}-x_0^n}{\Delta t_n}\right)
  \frac{u_1^{n+1}-u_{-1}^{n+1}}{2h_1^{n+1}} \nonumber \\
  & - \left(\frac{x_0^{n+1}-x_0^n}{\Delta t_n}\right)
  \frac{u_1^{n+1}-2u_0^{n+1}+u_{-1}^{n+1}}{(h_1^{n+1})^2} - s_0
  \sqrt{s_0} u_0^{n+1} = 0,
\end{align}
on the left boundary, and
\begin{align}\label{eq:LABC_right}
  & \frac{1}{\Delta t_n}\bigg(\frac{u_{N+1}^{n+1}-u_{N-1}^{n+1}}{2h_{N}^{n+1}} -
  \frac{u_{N+1}^{n}-u_{N-1}^{n}}{2h_{N}^{n}}\bigg) + 3 \sqrt{s_0}
  \frac{u_N^{n+1}-u_N^{n}}{\Delta t_n} + 3\left(s_0-\sqrt{s_0}\frac{x_N^{n+1}-x_N^n}{\Delta t_n}\right)
  \frac{u_{N+1}^{n+1}-u_{N-1}^{n+1}}{2h_N^{n+1}} \nonumber \\
  & - \left(\frac{x_N^{n+1}-x_N^n}{\Delta t_n}\right)
  \frac{u_{N+1}^{n+1}-2u_N^{n+1}+u_{N-1}^{n+1}}{(h_N^{n+1})^2} + s_0
  \sqrt{s_0} u_0^{n+1} = 0.
\end{align}
on the right boundary. Here two ghost points $x_{-1}$ and $x_{N+1}$
are used. On the other hand, if the observed domain is big enough or
else, Dirichlet boundary conditions are employed.

\subsection{Full algorithm}
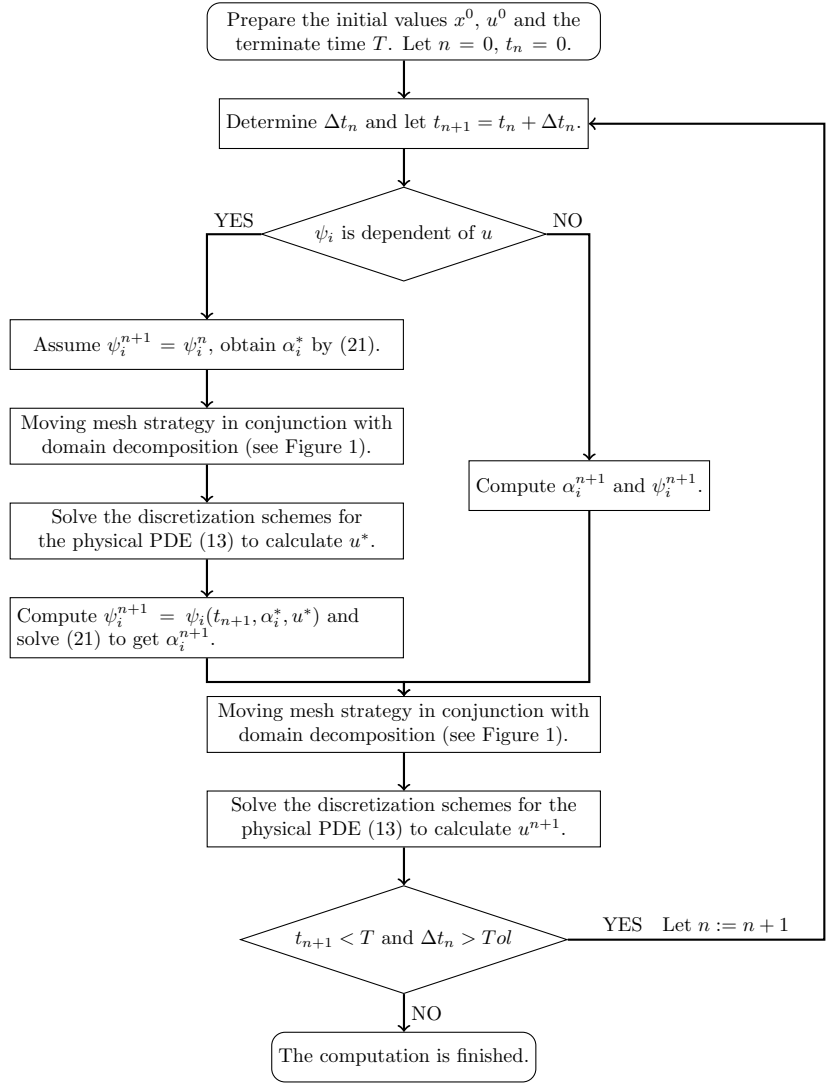
\begin{figure}[!htb]
  \centering
  \small
%\begin{minipage}{0.75\textwidth}
  \begin{tikzpicture}[scale = 0.75]
    \matrix[row sep=5mm,column sep=5mm, top color=white, bottom color=white, scale=0.75]{
      \node[startstop, text width=0.4\textwidth, scale=0.75, text centered] (start) {Prepare the initial values $x^0$, $u^0$ and the terminate time $T$. Let $n=0$, $t_n=0$.};\\
      \node[passprocess, scale=0.75] (timestep) {Determine $\Delta t_n$ and let $t_{n+1} = t_n + \Delta t_n$.};\\
      \node[decision, scale=0.75] (branchsource) {$\psi_i$ is dependent of $u$};\\
      \node[passprocess, text width=0.4\textwidth, left, scale=0.75, text centered] (pre_step1) {Assume $\psi_i^{n+1} = \psi_i^n$, obtain $\alpha_i^*$ by \eqref{eq:source_CN}.};\\
      \node[passprocess, left, text width=0.4\textwidth, scale=0.75, text centered] (preStep2) {Moving mesh strategy in conjunction with domain decomposition (see Figure \ref{fig:mesh_strategy}).};\\
      \node[passprocess, left, text width=0.4\textwidth, scale=0.75, text centered] (preStep3) {Solve the discretization schemes for the physical PDE \eqref{eq:heat_logical} to calculate $u^*$.};\\
      \node[passprocess, text width=0.4\textwidth, left, scale=0.75] (corstep1) {Compute $\psi_i^{n+1}=\psi_i(t_{n+1}, \alpha_i^*, u^*$) and solve \eqref{eq:source_CN} to get $\alpha_i^{n+1}$.};\\
      \node[passprocess, text width=0.4\textwidth, scale=0.75, text centered] (corstep2) {Moving mesh strategy in conjunction with domain decomposition (see Figure \ref{fig:mesh_strategy}).};\\
      \node[passprocess, text width=0.4\textwidth, scale=0.75, text centered] (corstep3) {Solve the discretization schemes for the physical PDE \eqref{eq:heat_logical} to calculate $u^{n+1}$.};\\
%      \node[] {};\\
      \node[decision, scale=0.75, text centered] (branchtime) {$t_{n+1}<T$ and $\Delta t_n>Tol$};\\
      \node[startstop, scale=0.75] (end) {The computation is finished.}; \\
    };
    \node[passprocess, scale=0.75] (pdeStep1) at ($(preStep3.north east)+(3.25cm,3mm)$) {Compute $\alpha_i^{n+1}$ and $\psi_i^{n+1}$.};
  \begin{scope}[every path/.style=line]
    \path (start) -- (timestep);
    \path (timestep) -- (branchsource);
    \path (branchsource) -| node[near start, above, scale=0.75] {YES} (pre_step1);
    \path (branchsource) -| node[near start, above, scale=0.75] {NO} (pdeStep1);
    \path (pre_step1) -- (preStep2);
    \path (preStep2) -- (preStep3);
    \path (preStep3) -- (corstep1);
    \path (corstep2) -- (corstep3);
    \path (corstep1) |- ($(corstep2.north)+(0,2.5mm)$) -- (corstep2);
    \path (corstep3) -- (branchtime);
    \path (pdeStep1) |- ($(corstep2.north)+(0,2.5mm)$) -- (corstep2);
    \path (branchtime) -- node[midway, right, scale=0.75] {NO} (end);
    \path (branchtime) -- node[midway, above, scale=0.75] {YES\quad Let $n:=n+1$} ($(branchtime.east)+(4.5,0)$) |- (timestep);
  \end{scope}
  \end{tikzpicture}
%\end{minipage}
  \caption{Full algorithm for numerical solution of the model problem \eqref{eq:heat_eq}-\eqref{eq:source_eq}.}
%  \caption{Full algorithm for the model problem \eqref{eq:heat_eq}$-$\eqref{eq:bound_cond} and \eqref{eq:source_eq}.}
  \label{fig:full_algorithm}
\end{figure}

We close this section with a full algorithm in Figure
\ref{fig:full_algorithm} for the model problem
\eqref{eq:heat_eq}$-$\eqref{eq:bound_cond} and
\eqref{eq:source_eq}. Here the choice of the time step $\Delta t_n$
will be determined in the following concrete examples, and $Tol>0$ is
set to be $10^{-16}$.

\section{Numerical examples}\label{sec:example}
In this section, we present some numerical examples to verify the
convergence rate and illustrate efficiency of the full algorithm in
Figure \ref{fig:full_algorithm}.
% The first one is a constructed example for the nonlinear moving
% interface problem \cite{li1997immersed,beyer1992analysis}. The rest
% are from traveling heat source problem with the solution may be
% blow-up \cite{zhu2010numerical,ma2009moving,hu2011moving,hu2012zh}.

\begin{example}
  We consider a nonlinear moving interface problem with the following
  exact solution
\begin{align}
  \label{eq:ex_true_u}
  u(x,t)=
  \begin{cases}
    \sin( \omega_1 x ) e^{- \omega_1^2 t}, \quad & x\leq \alpha_0(t), \\
    \sin( \omega_2 (1-x)) e^{- \omega_2^2 t}, \quad & x\geq \alpha_0(t),
  \end{cases}
\end{align}
for some choice of $\omega_1$ and $\omega_2$. The interface
$\alpha_0(t)$ is determined by solving the scalar equation
\begin{equation}
  \label{ex:interface_eq}
  \sin(\omega_1 \alpha_0 ) e^{-\omega_1^2t} = \sin(\omega_2(1-\alpha_0)) e^{-\omega_2^2 t},
\end{equation}
so that $u(x,t)$ is continuous across the interface.
\end{example}

The equation \eqref{ex:interface_eq} has a unique solution on $[0,1]$
if we take, for example, $\pi<\omega_1,~\omega_2<2\pi$. Then we have
the ordinary differential equation for the motion of the interface
\begin{align}
  \label{ex:interface_ode}
  \frac{\mathrm{d} \alpha_0}{\mathrm{d}t} = \frac{(\omega_1^2-\omega_2^2) u(\alpha_0, t)}{u_x(\alpha^-_0, t) - u_x(\alpha_0^+, t)}.
\end{align}
Based on the jump conditions, the source function $F_0(t,x,u)$ is
\begin{align}
  \label{ex:source_fun}
  F_0(t,x,u) &= -[u_x]_{\alpha_0} = u_x(\alpha^-_0, t) - u_x(\alpha_0^+, t)\nonumber\\
  & = \omega_1\cos(\omega_1\alpha_0)e^{-\omega_1^2t} + \omega_2 \cos(\omega_2(1-\alpha_0))e^{-\omega_2^2t}.
\end{align}
Same as in \cite{beyer1992analysis}, we take $\omega_1=5\pi/4$,
$\omega_2=7\pi/4$. The observed domain is set by $[0,1]$, where the
initial position of the interface is $\alpha_0(0)=0.58333$. Since we
have the exact solution, Dirichlet boundary conditions are
employed. In this example, we simply use the uniform time step,
i.e. $\Delta t_n \equiv const$, and the total number of the time
meshes is $L$. The monitor function for MMPDE6 \eqref{eq:MMPDE6} takes
the form
\begin{equation}\label{ex:monitor}
  M(x,t) = (1-\theta) \bigg|\frac{\partial u}{\partial
    x}\bigg| + \theta((x-\alpha_0(t))^2+\varepsilon)^{-1/4},
\end{equation}
where $0<\theta<1$, $0<\varepsilon\ll 1$. This is consistent with the
choice in \cite{hu2011moving,hu2012zh}. In practice, smoothing the
monitor function can improve the accuracy of the numerical solution,
and we utilize the smoothing technique proposed in
\cite{huang1994moving}. Here the parameters in MMPDE6
\eqref{eq:MMPDE6} and the monitor function \eqref{ex:monitor} are
given by $\tau=10^{-3}$, $\theta=0.5$, and $\varepsilon=10^3/N^4$.

Since backward Euler scheme is used to solve the physical PDE in this
paper, the truncation error for time discretization is only
first-order. To verify our algorithm has a second-order convergence
rate for space, the number of $L$ should be fourfold when $N$ is
double in the convergence test.
% In the test, we use the solution of a zero-finding MATLAB function
% \emph{fzero} for \eqref{ex:interface_eq} as the exact position of
% the interface. For the true solution is known, the speed of the
% interface can also be evaluated exactly. We compare the error of the
% solution and the interface in the time $T=0.1$, for the full
% predictor-corrector algorithm and the algorithm where $\alpha_0(t)$
% and $\alpha_0'(t)$ are exactly calculated.
Computational results with different number of $N$ and $L$ at the time
$T=0.1$ are listed in Table \ref{table:err_result}, where the errors
are defined as
\begin{align*}
  E_{N,L} = \parallel U - u_e \parallel_\infty,
  \quad \tilde{E}_{N,L} = \parallel \tilde{U} - u_e \parallel_\infty,
  \quad E_{N,L}^\alpha = \mid \tilde{\alpha}_0 - \alpha_0^* \mid.
\end{align*}
Here, $u_e$ is the true solution, $\alpha_0^*$ used as the exact
interface is the solution of a zero-finding MATLAB function
\emph{fzero} for \eqref{ex:interface_eq}. The numerical solution $U$
is obtained by the algorithm where $\alpha_0(t)$ and $\alpha_0'(t)$
are exactly calculated. And $\tilde{U}$, $\tilde{\alpha}_0$ represent
respectively the numerical solution and interface, solved with the
full predictor-corrector algorithm. The ratios in Table
\ref{table:err_result} are $E_{2N,4L}/E_{N,L}$,
$\tilde{E}_{2N,4L}/\tilde{E}_{N,L}$ and
$E_{2N,4L}^\alpha/E_{N,L}^\alpha$, respectively. It is shown that our
algorithm solves the solution and the interface very well, and has a
second-order convergence rate for space,
i.e. $O(1/N^2)$. Additionally, compared the corresponding results in
\cite{li1997immersed}, our algorithm is better than the method
proposed in \cite{li1997immersed}.

\begin{table}[!htb]
  \caption{Error and convergence rates at $T=0.1$.}% 数据来源文件："images\ode_source1_MovingMesh"
  \label{table:err_result}%这个标签放到 \begin{table} 这一行引用竟然序号不对！什么缘故？
  \centering
  \tiny
  \colorbox{white}{
    \begin{tabular}{lllllll}
      \hline
      N,~L          & $E_{N,L}$   & ratio       & $\tilde{E}_{N,L}$   & ratio       & $E_{N,L}^\alpha$       & ratio \\
      \hline
      40, ~40       & 1.0931e-2    & -          & 1.3721e-2    & -          & 8.2720e-3    & - \\
      80, ~160      & 2.6996e-3    & 0.24697    & 3.3908e-3    & 0.24713    & 2.0540e-3    & 0.24830 \\
      160, ~640     & 6.6945e-4    & 0.24798    & 8.4144e-4    & 0.24816    & 5.1038e-4    & 0.24848 \\
      320, ~2560    & 1.6687e-4    & 0.24927    & 2.0981e-4    & 0.24935    & 1.2732e-4    & 0.24947 \\
      640, ~10240   & 4.1678e-5    & 0.24976    & 5.2408e-5    & 0.24978    & 3.1807e-5    & 0.24982 \\
      1280, ~40960  & 1.0416e-5    & 0.24992    & 1.3098e-5    & 0.24993    & 7.9500e-6    & 0.24994 \\
      2560, ~163840 & 2.6038e-6    & 0.24998    & 3.2743e-6    & 0.24998    & 1.9874e-6    & 0.24998 \\
      \hline
    \end{tabular}}
\end{table}

\begin{figure}[!htb]
  \centering
  {\includegraphics[width=0.335\textwidth]{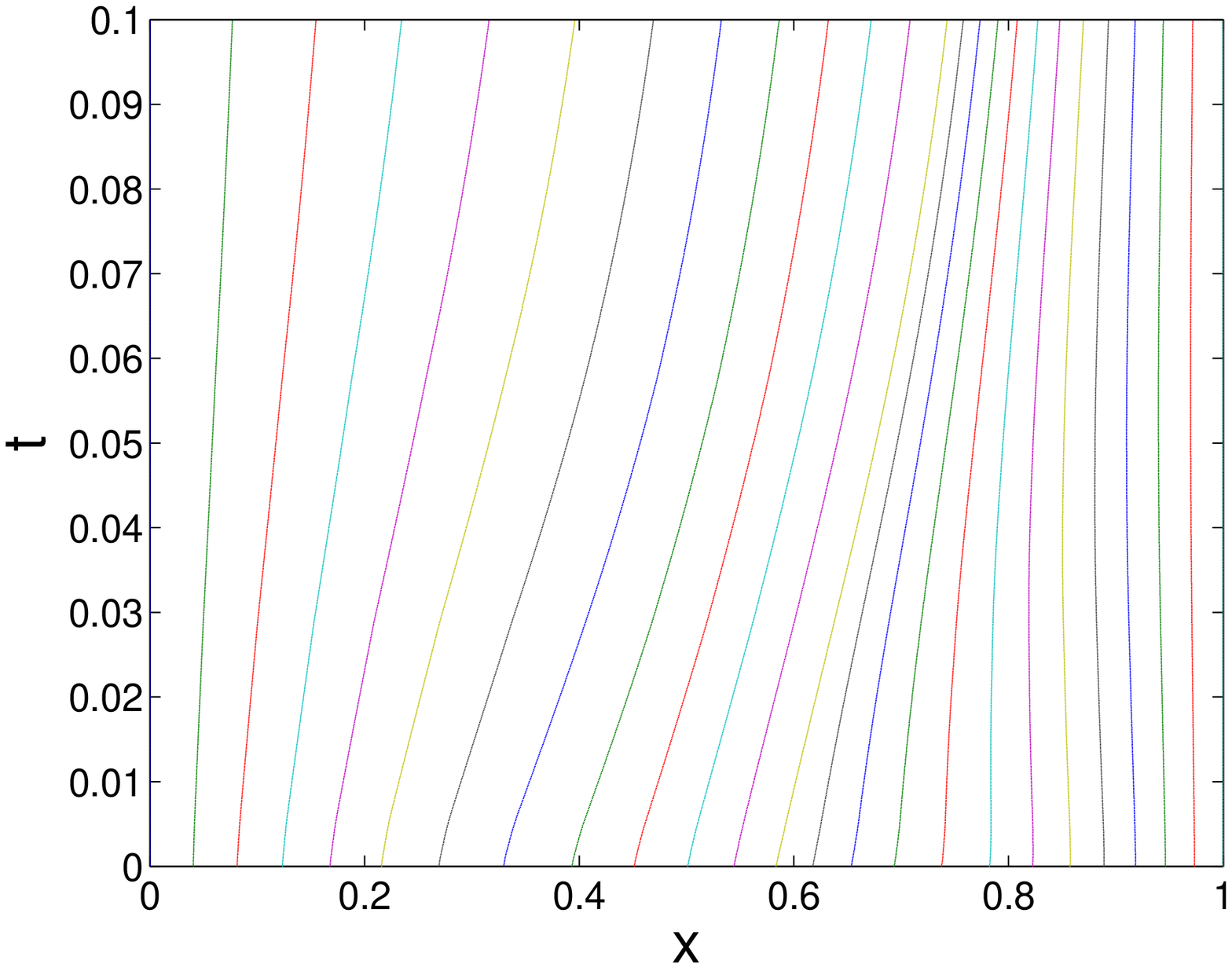}}%\hfill
  {\includegraphics[width=0.33\textwidth]{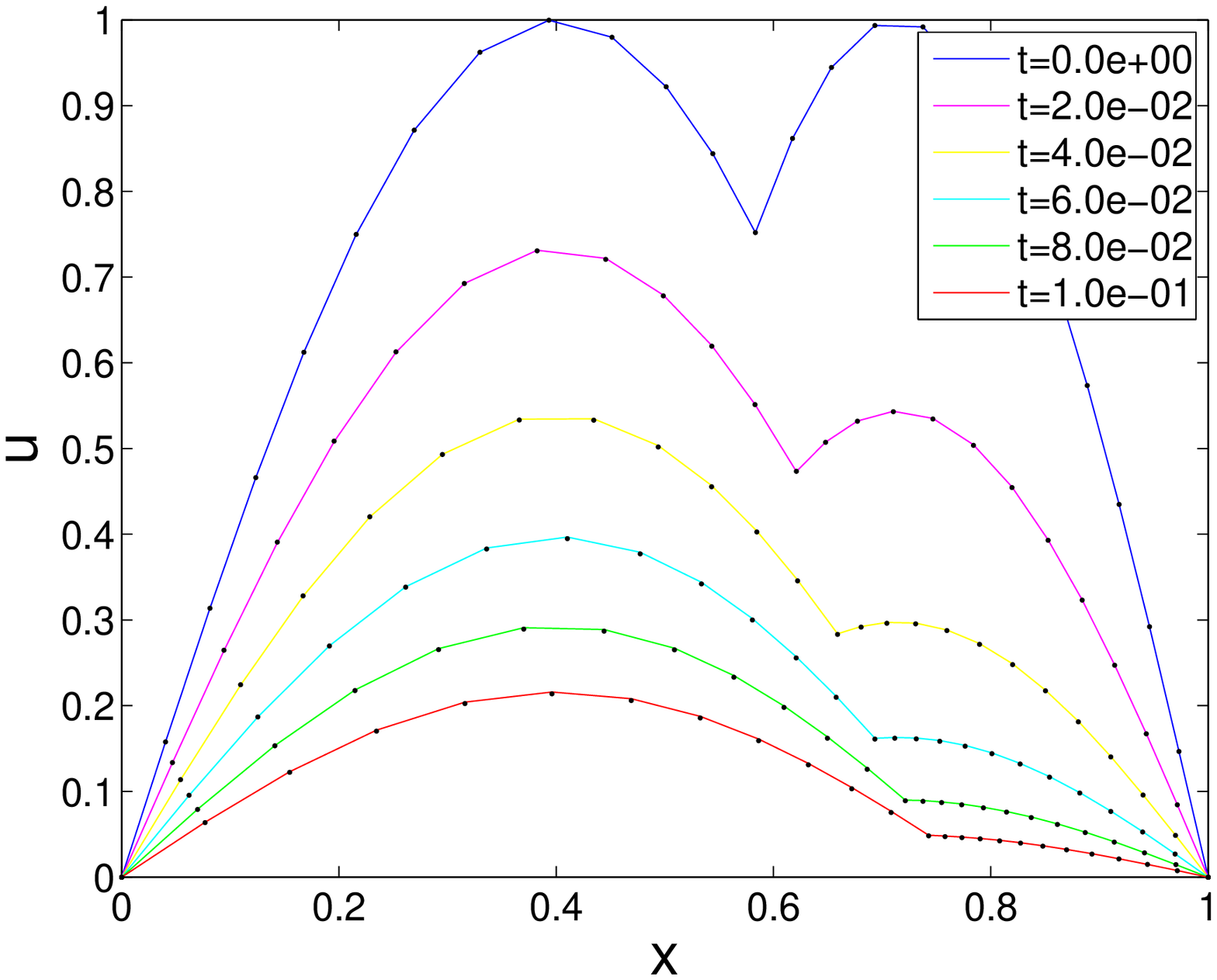}}%\\
  {\includegraphics[width=0.33\textwidth]{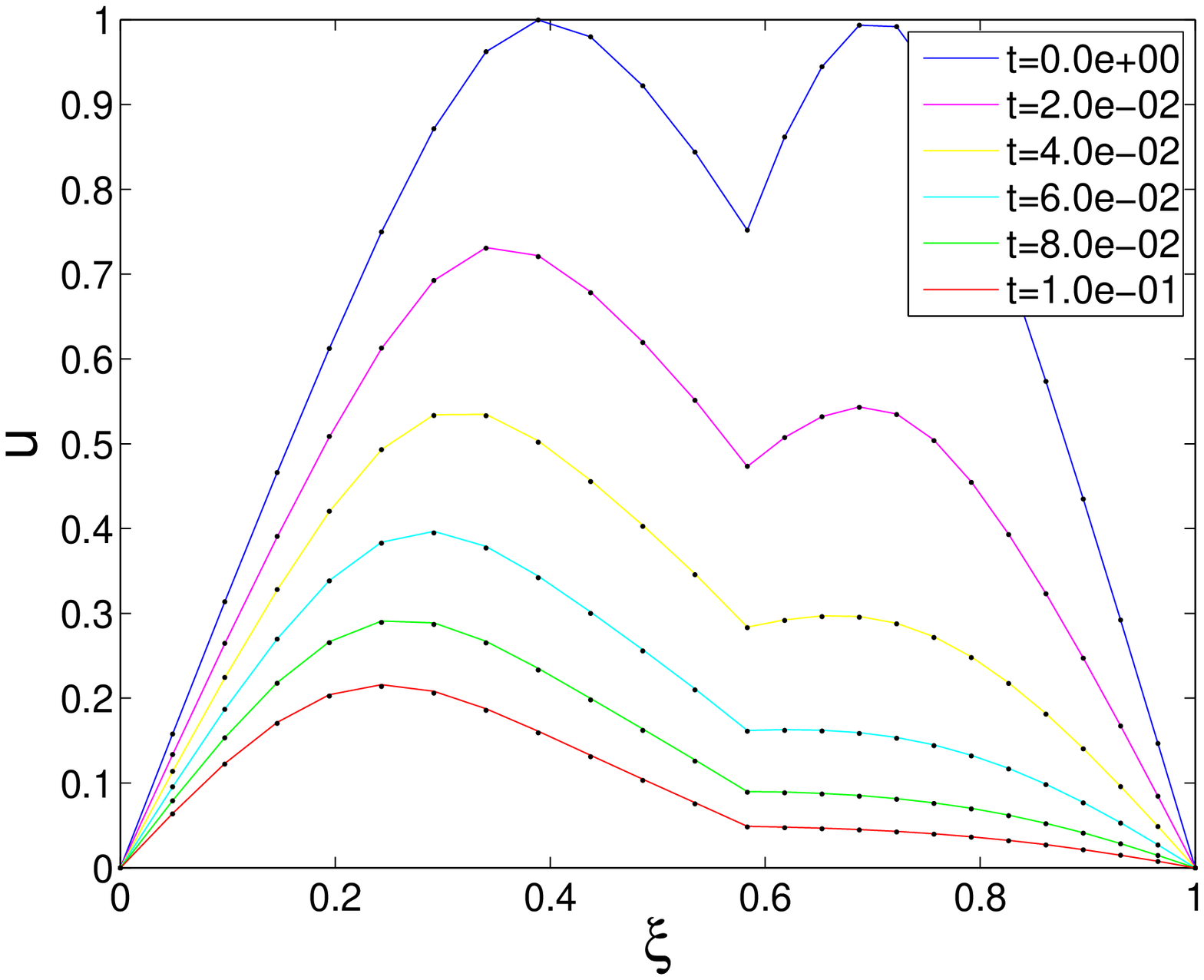}}%\hfill
  \caption{\small Mesh trajectories and the profiles of $u$ in
    physical variable and computational variable (from left to right)
    as time changes with $N=24$. The solid lines are the computed
    solution and the dots are the exact solution on the mesh.}
  \label{fig:ex1}
\end{figure}

Figure \ref{fig:ex1} presents the profiles of the solution in physical variable and computational variable and the evolving mesh from $t=0$ to $t=0.1$. The number of the mesh is $N=24$, with half mesh points on each side of the interface. We can see that we get excellent resolution of the example even with a grid as coarse as $N=24$.

\bigskip
The rest examples are from traveling heat sources problems with the solution may be blow-up \cite{zhu2010numerical,ma2009moving,hu2011moving,hu2012zh}. If not specifically pointed out, the initial value is given by
\begin{align}\label{eq:init_value}
  u(x,0) = \left\{ \begin{array}{ll}
      \cos^2(\pi x/2), & -1<x<1,\\ 0, &
      \text{otherwise},
    \end{array} \right.
\end{align}
the observed domain is set by $[-10,10]$ with $u(-10,t)=u(10,t)=0$, and the source functions $F_i(t,x,u)$ are simply specified by
\begin{equation}\label{eq:F_source}
  F_i(t,x,u) = 1 + u^2, \qquad i=0,1,\ldots,q-1.
\end{equation}
The resulting nonlinear system is solved by Newton iteration with the tolerance $tol = 10^{-8}$.

The monitor function for MMPDE6 \eqref{eq:MMPDE6} takes the form
\begin{equation}\label{eq:monitor}
  M(x,t) = \theta_{q+1} u^p + \theta_q \bigg|\frac{\partial u}{\partial
    x}\bigg| + \sum_{i=0}^{q-1}\theta_i((x-\alpha_i(t))^2+\varepsilon)^{-1/4},
\end{equation}
where the parameters $0<\theta_i<1$, $\sum_{i=0}^{q+1}\theta_i=1$, $0<\varepsilon\ll 1$, $p>0$ will be determined later. For non-blowup case, the following graded time steps
\cite{ma2009moving,zhu2010numerical}
\begin{equation*}
  t_n=\bigg(n\frac{T}{L}\bigg)^2, \quad n = 0, 1, \ldots, L,
\end{equation*}
are used with $[0,T]$ the time integration interval and $L$ the number of time meshes. While for blow-up case, the time step $\Delta t_n = t_{n+1}-t_{n}$ is chosen to be \cite{ma2009moving,budd2001scaling}
\begin{equation*}
  \Delta t_n =\min\left \{\mu, \frac{\mu}{\Big(\max_{j}\big\{u_j^n\big\}+\varepsilon\Big)^2}\right \},
\end{equation*}
where $\varepsilon$ is same in the monitor function, $\mu$ is a small positive constant with $\mu=10^{-3}$ in the test.

\begin{example}[Linear moving sources]
  We consider that all sources move with a constant velocity $k$, and
  the position of the $i$-th source has a constant distance $d_i$ to
  the $0$-th source, i.e.,
  \begin{equation*}
    \alpha_i'(t) = k, \quad \alpha_i(0)=d_i,\qquad i=0,1,\ldots,q-1,
  \end{equation*}
where $d_0=0$.
\end{example}

Since our method is trivial for multi-sources case, only $q=1,2$ are considered. Different velocities $k$ are investigated in \cite{zhu2010numerical,ma2009moving,hu2011moving} and we specify $k=2$ here. With this velocity, blow-up would occur for $q=2$, and be avoided for $q=1$. The parameters are set by $\tau=10^{-3}$, $\theta_0=0.9$, $\theta_1=0.1$, and $\varepsilon=10^3/N^4$ for $q=1$, while $\tau=5\times 10^{-4}$, $\theta_0=\theta_1=0.3$, $\theta_3=0.4$, $p=2$, $\varepsilon = 10^{-5}$, and $d_1=2.5$ for $q=2$, respectively.

\begin{figure}[!htb]
  \centering
  {\includegraphics[width=0.33\textwidth]{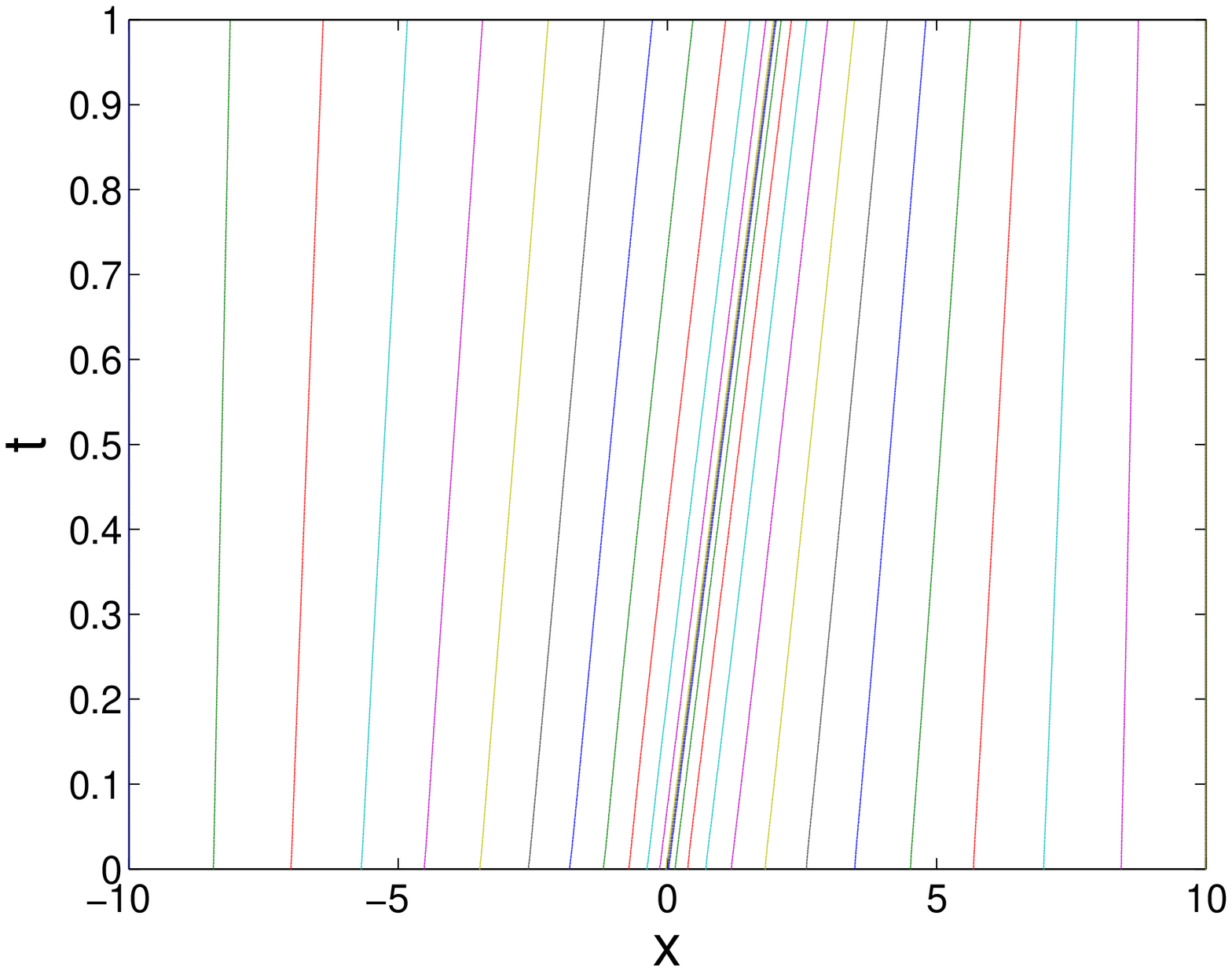}}%\hfill
  {\includegraphics[width=0.33\textwidth]{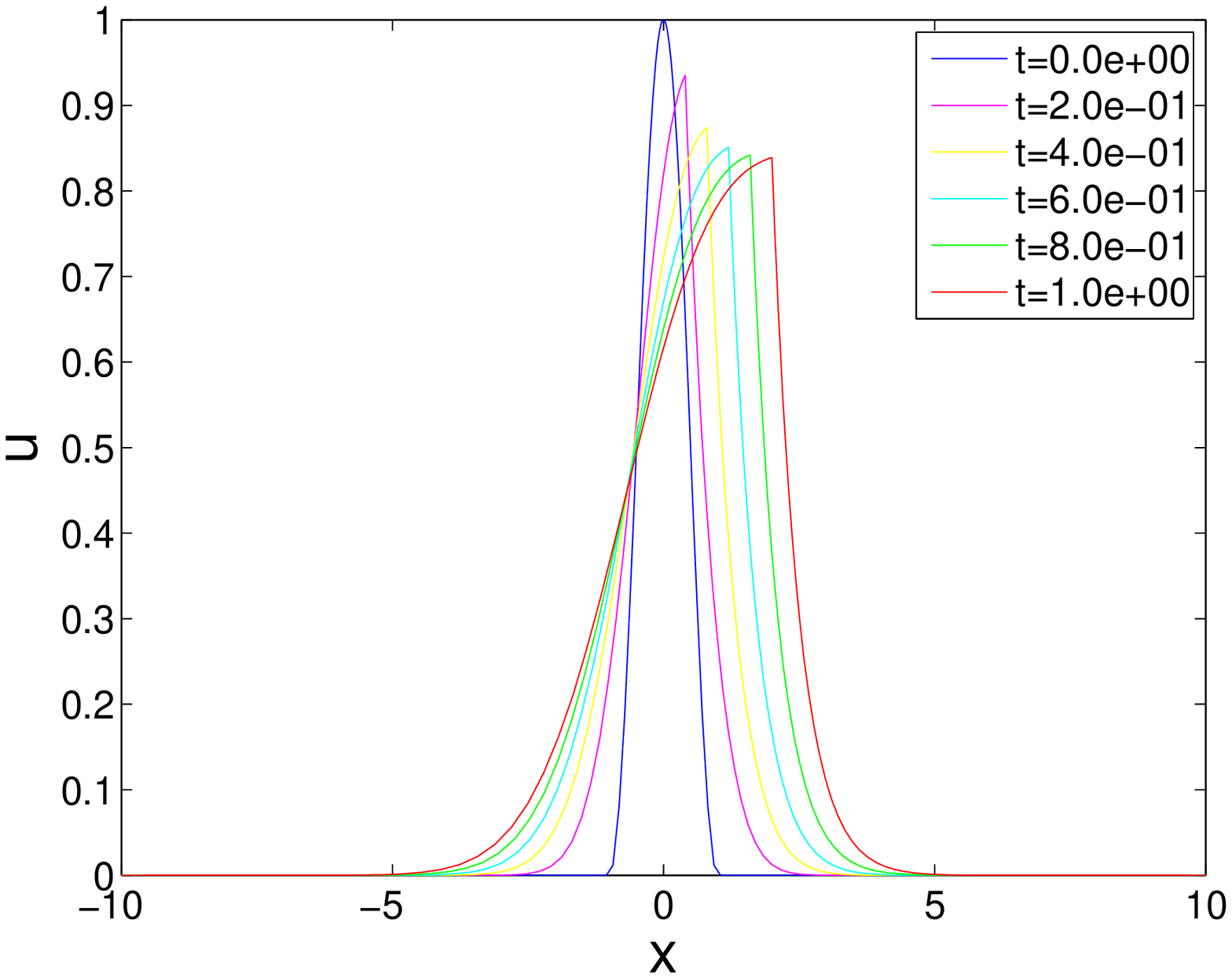}}%\\
  {\includegraphics[width=0.33\textwidth]{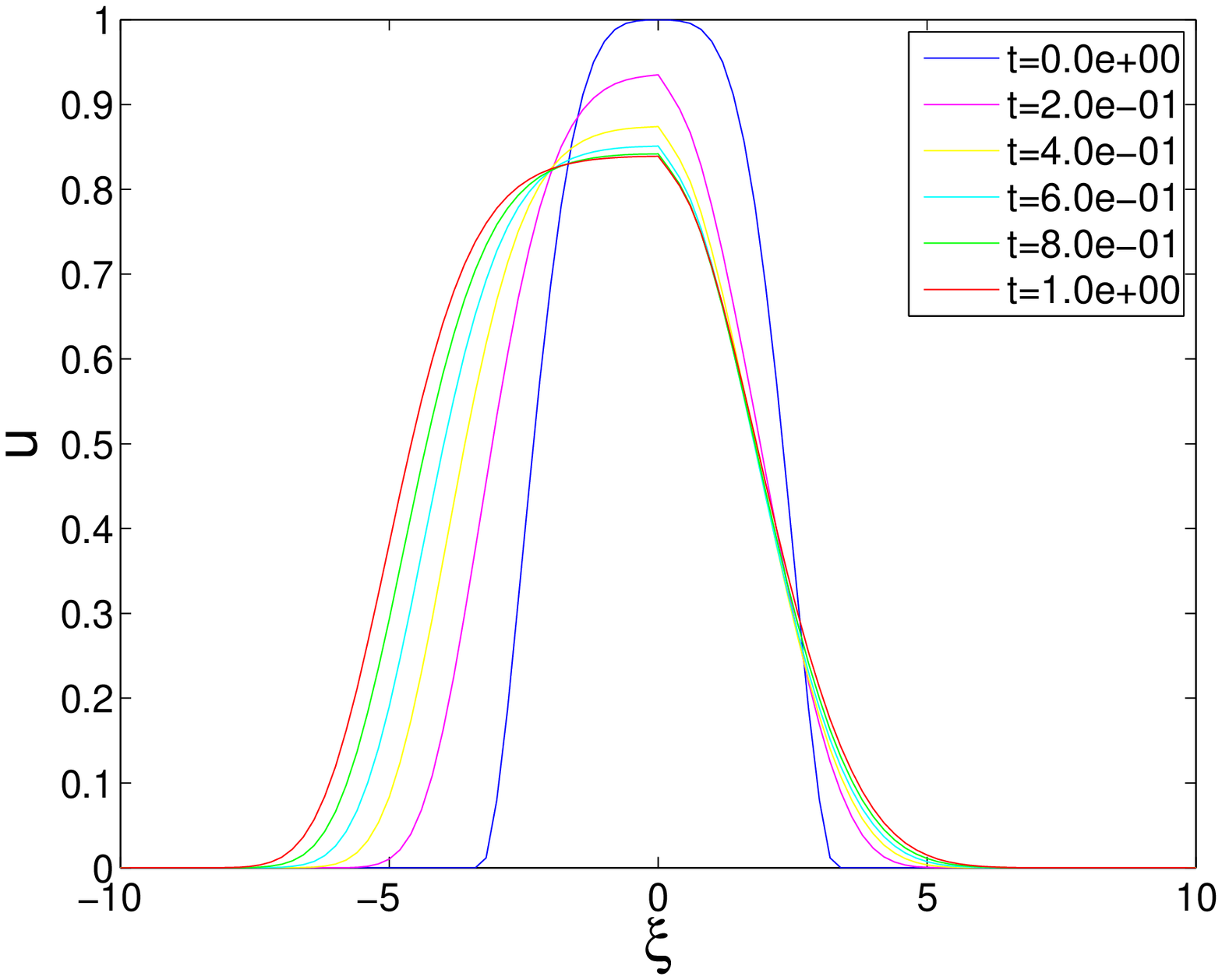}}%\hfill
  \caption{\small Mesh trajectories and the profiles of $u$ from $t=0$
    to $t=1.0$ for one source case with $N=100$.}
  \label{fig:linear_nonblowup}
\end{figure}

\begin{figure}[!htb]
  % 就一张网格的图上 1 M, 结果生成的 pdf 还是有 12 M... 没办法, 还是把
  % 网格图用 png 转换过来, 这样只好重新调整原来设好的比例了, 这个时
  % 候 pdf 还有 1.5 M 左右.
  \centering
  {\includegraphics[width=0.36\textwidth]{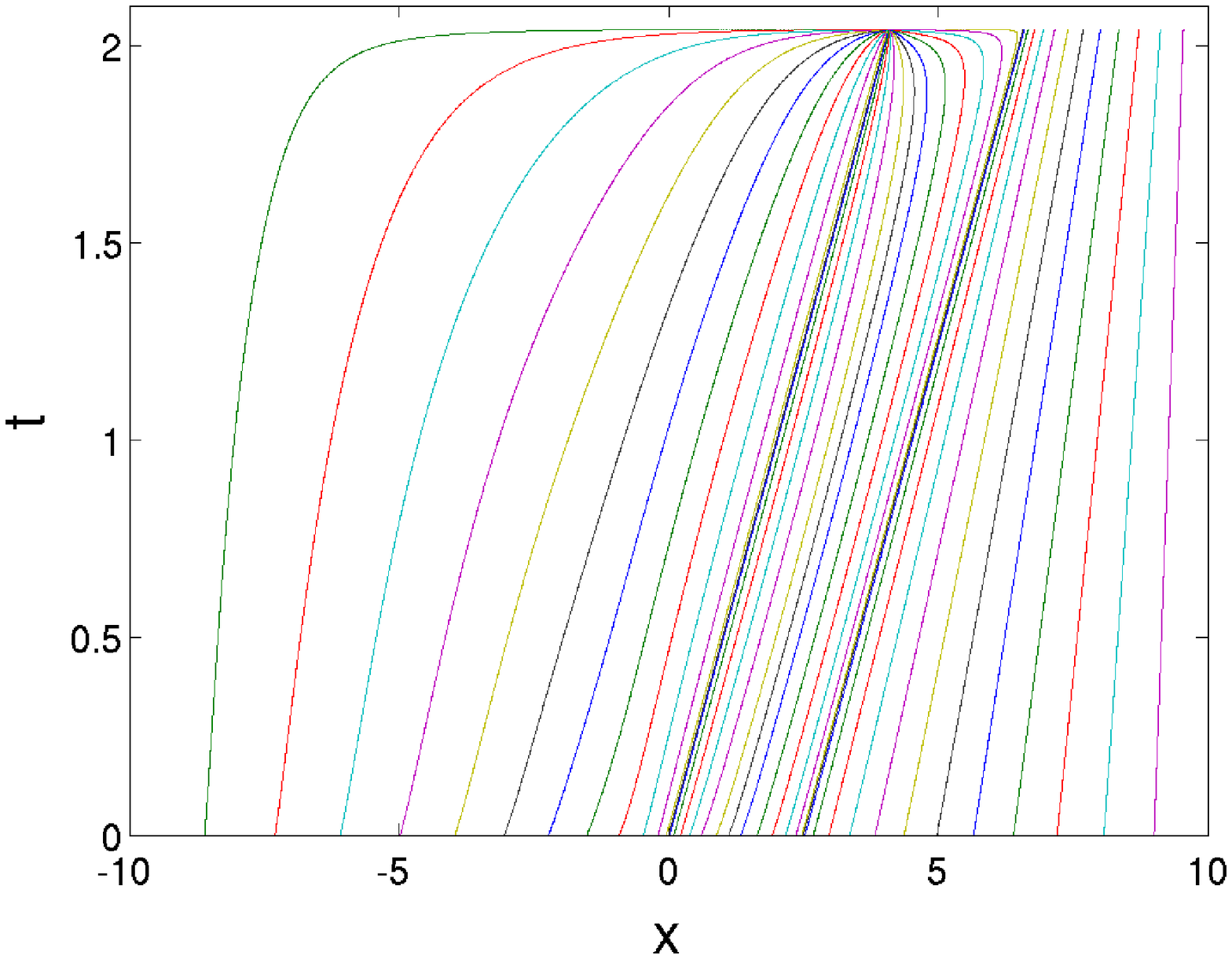}}%\hfill
  {\includegraphics[width=0.32\textwidth]{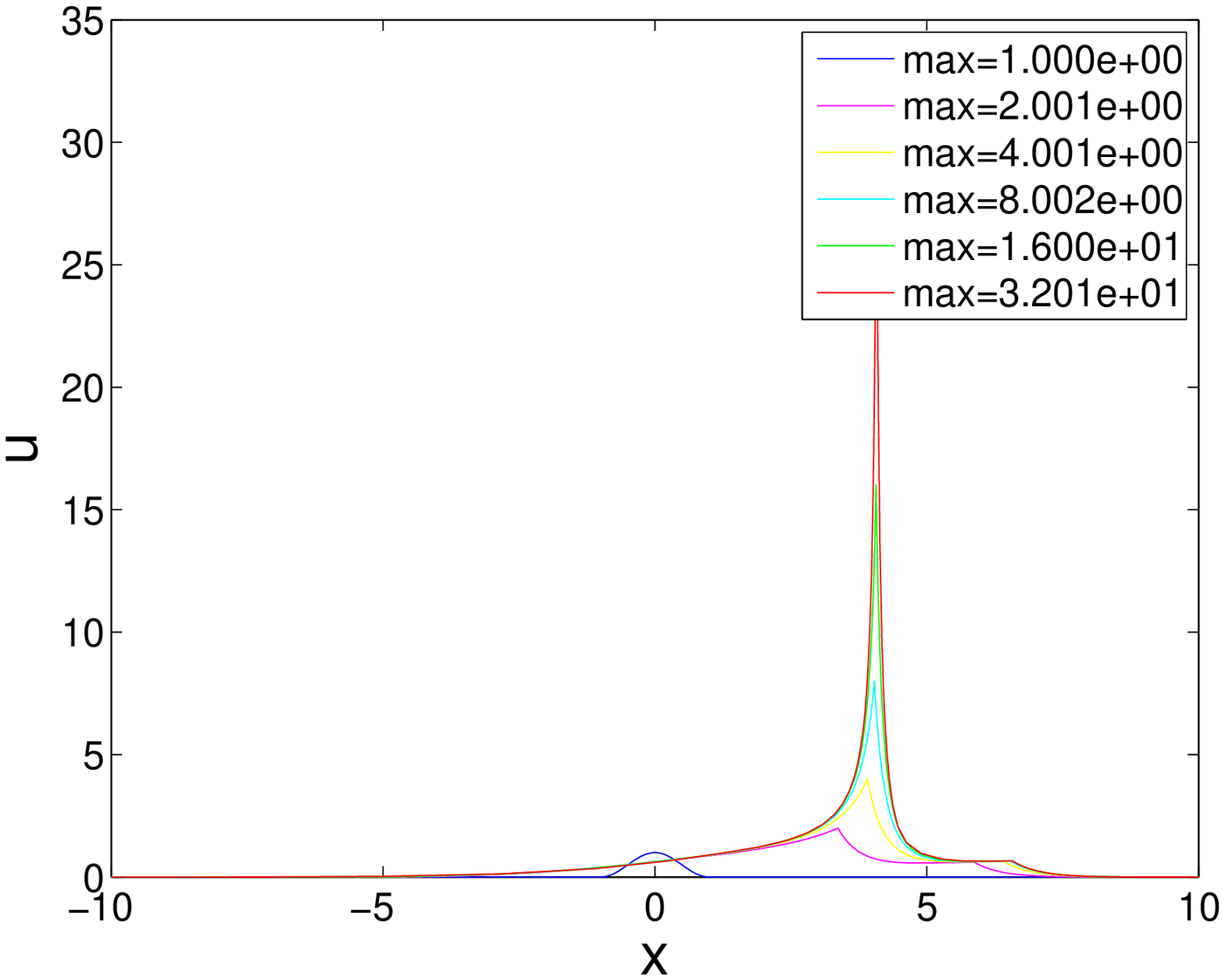}}%\\
  {\includegraphics[width=0.32\textwidth]{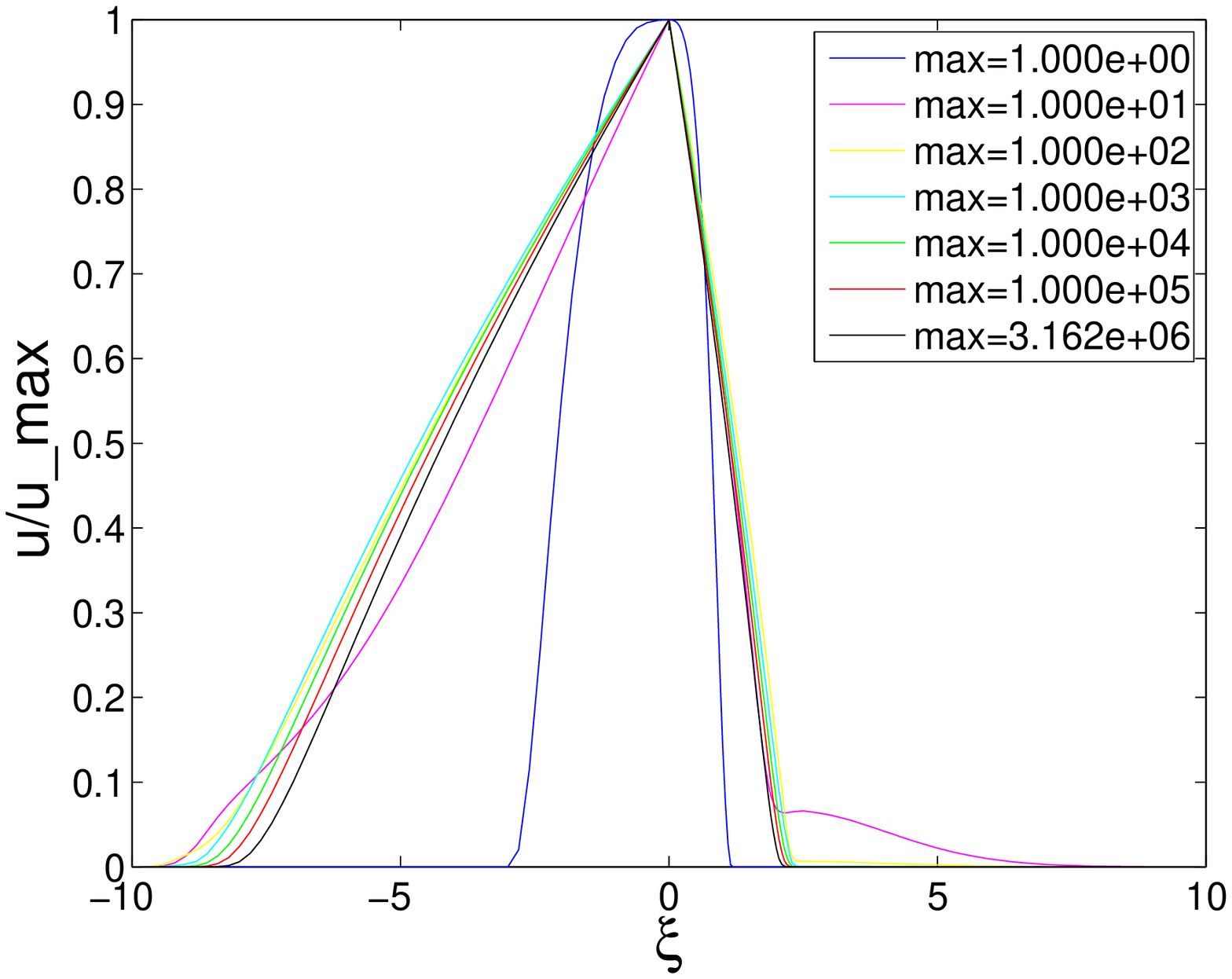}}%\hfill
  \caption{\small Mesh trajectories and the profiles of $u$ for two sources case with $N=150$.}
  \label{fig:linear_blowup}
\end{figure}

The profiles of the computed solution in physical variable and computational variable and the evolving mesh are presented in Figure \ref{fig:linear_nonblowup} for $q=1$ and in Figure \ref{fig:linear_blowup} for $q=2$. For simplicity, each subdomain has $50$ mesh points, i.e. $N=100$ for $q=1$ and $N=150$ for $q=2$. The numerical results are coincide with that in \cite{zhu2010numerical,ma2009moving}, and the blow-up time is $2.039708648680643$ at the first source $x=4.079417297361286$, corresponding to the maximum value of $u_{\max}=3.16\times 10^6$.

\begin{example}[Sin-type moving sources]
  We now consider two sources case, in which the sources move
  periodically with the same speed while separated by a constant
  distance $d_1=2.5$, that is,
  \begin{equation*}
    \alpha_0'(t) = \alpha_1'(t) = A \cos(\pi t),\quad \alpha_0(0) = 0.
  \end{equation*}
\end{example}

\begin{figure}[!htb]
  \centering
  {\includegraphics[width=0.36\textwidth]{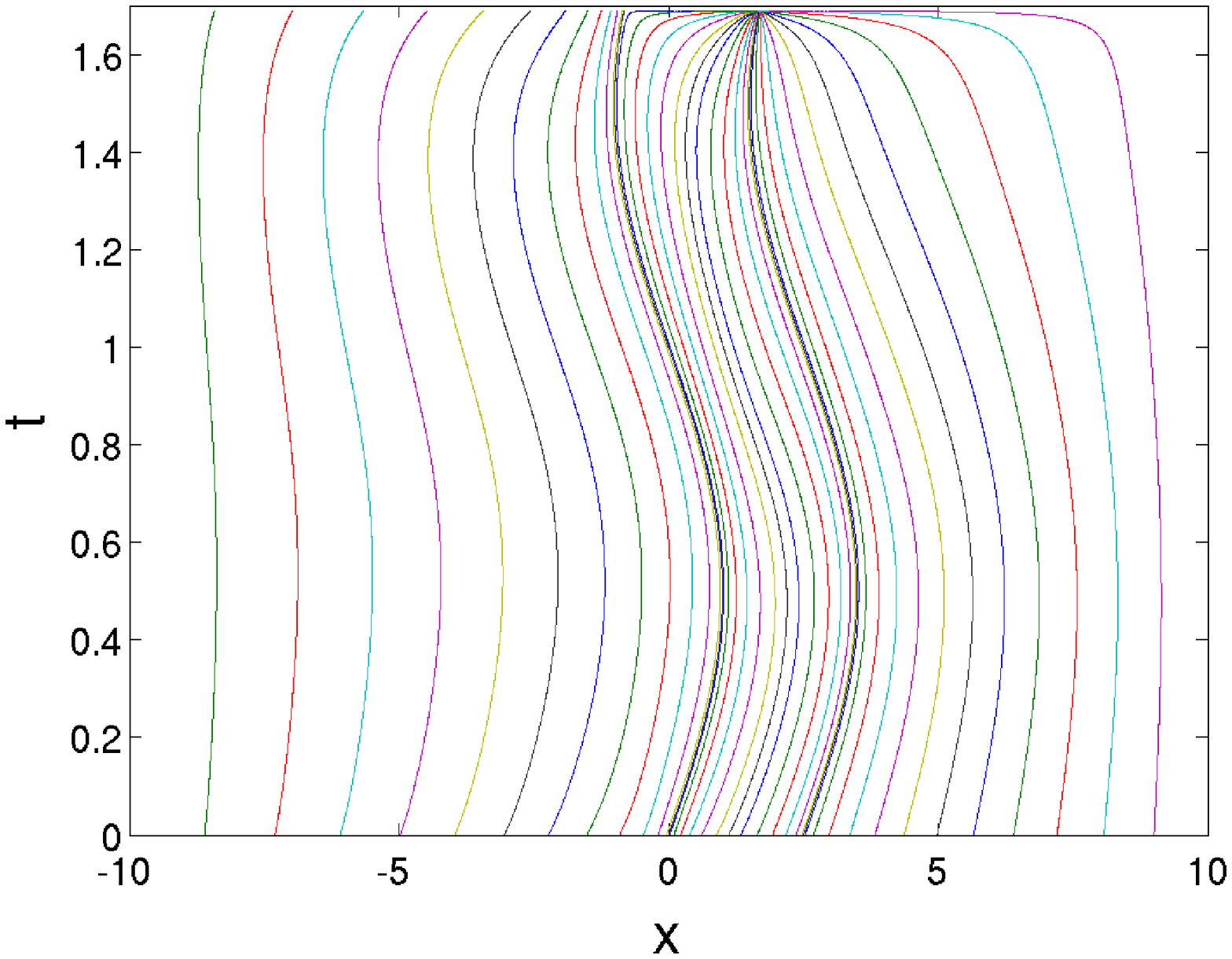}}%\hfill
  {\includegraphics[width=0.32\textwidth]{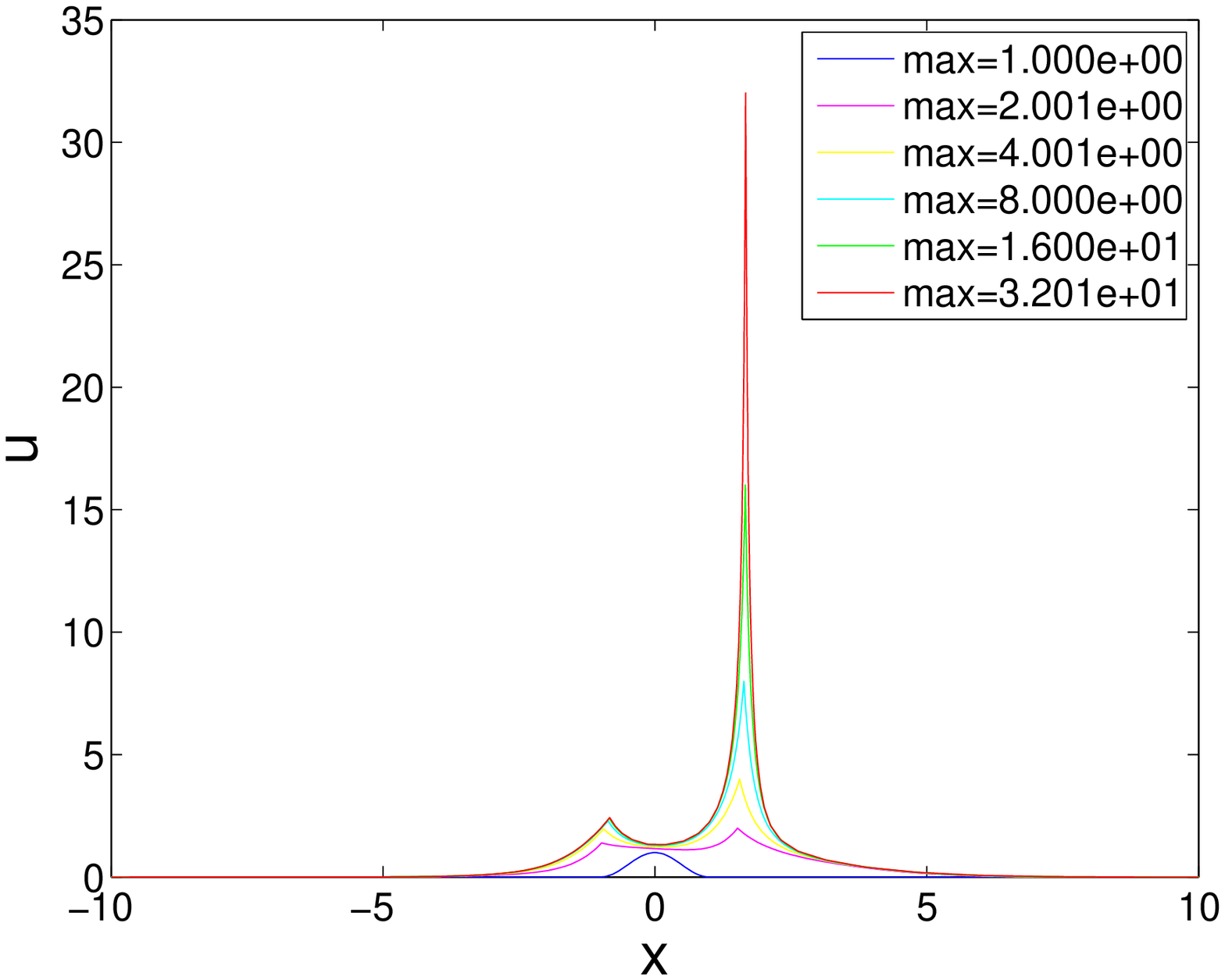}}%\\
  {\includegraphics[width=0.32\textwidth]{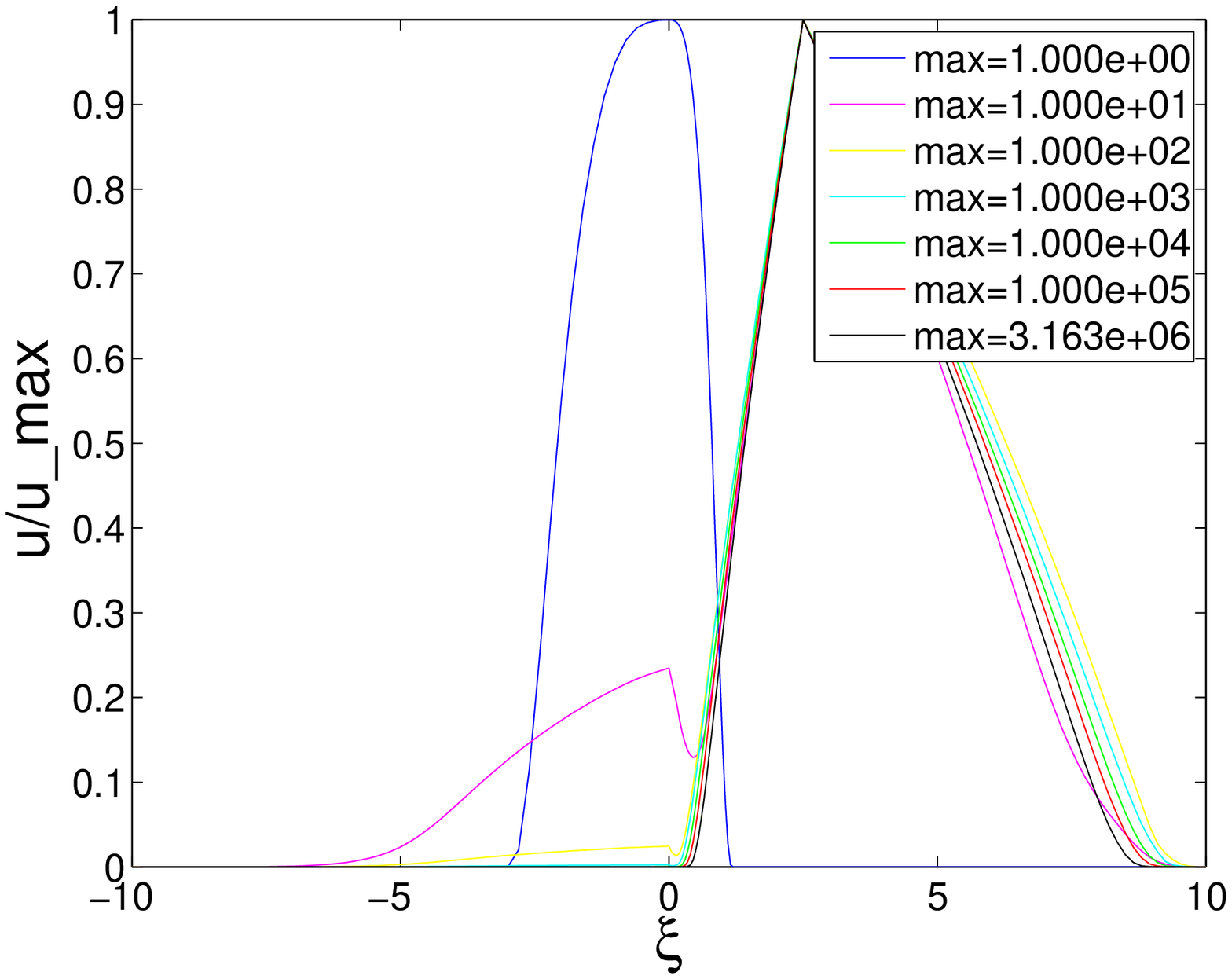}}%\hfill
  \caption{\small Mesh trajectories and the profiles of $u$ for $q=2$, $A=\pi$ with $N=150$.}
  \label{fig:sin1_source2}
\end{figure}

The Blow-up phenomenon is studied in \cite{hu2011moving} for different amplitudes $A$. Here we only give the numerical results for $A=\pi$ (see Figure \ref{fig:sin1_source2}), since all results are similar to those in \cite{hu2011moving}. All parameters are chosen the same as those in last example. The blow-up occurs at $t=1.689611393639939$ on the second source with the maximum value of $u_{\max}=3.16\times 10^6$.

\begin{example}[Symmetric periodic moving sources]
  We consider the case for two sources, which move periodically and
  symmetrically. The motion are described by
  \begin{equation*}
    \alpha_0'(t) = A \cos(\pi t),\quad \alpha_0(0) = -2.0,
  \end{equation*}
and $\alpha_1(t)=-\alpha_0(t)$, with e.g. $A=\pi$.
\end{example}

\begin{figure}[!t]%[!htb]
  \centering
  {\includegraphics[width=0.36\textwidth]{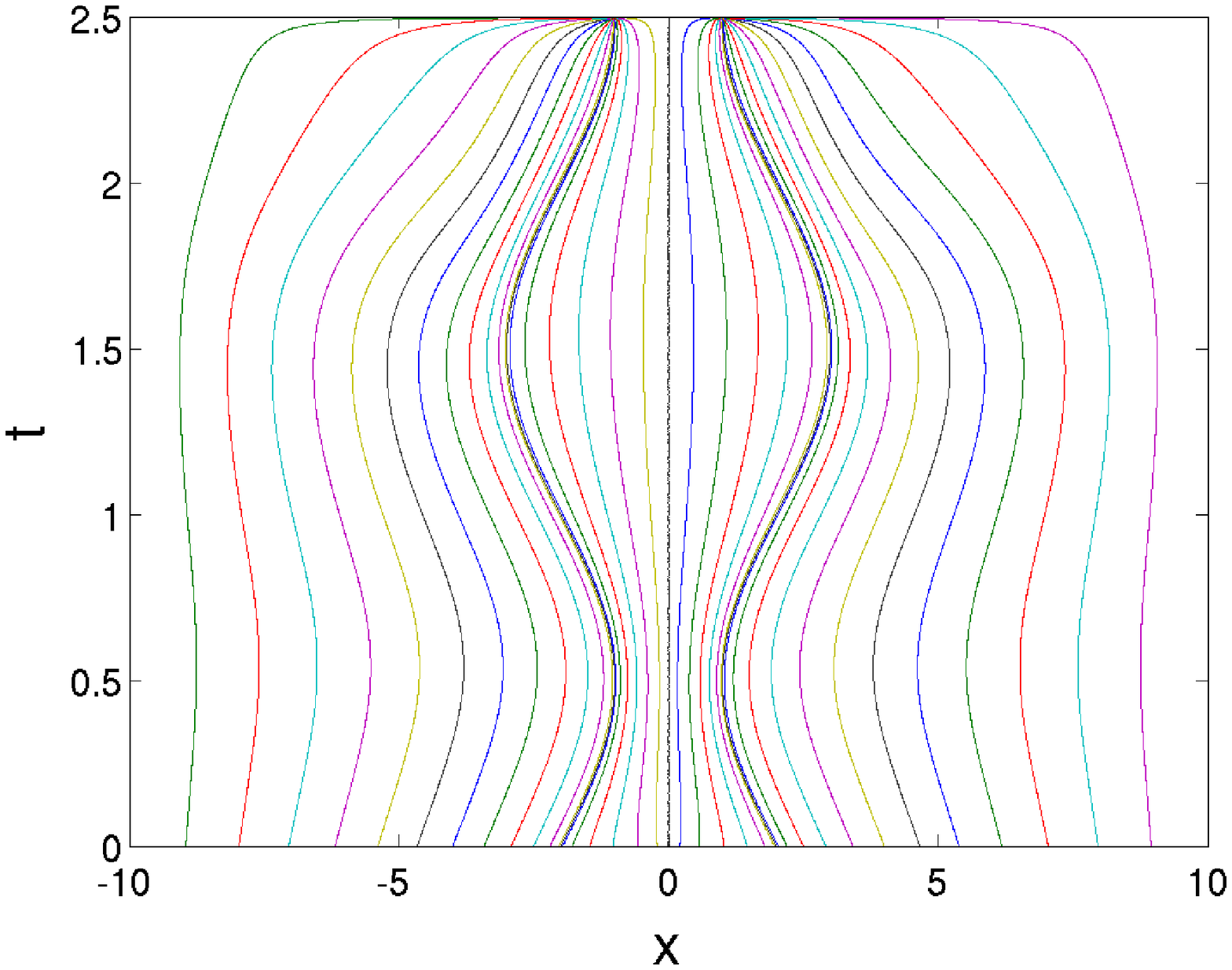}}%\hfill
  {\includegraphics[width=0.32\textwidth]{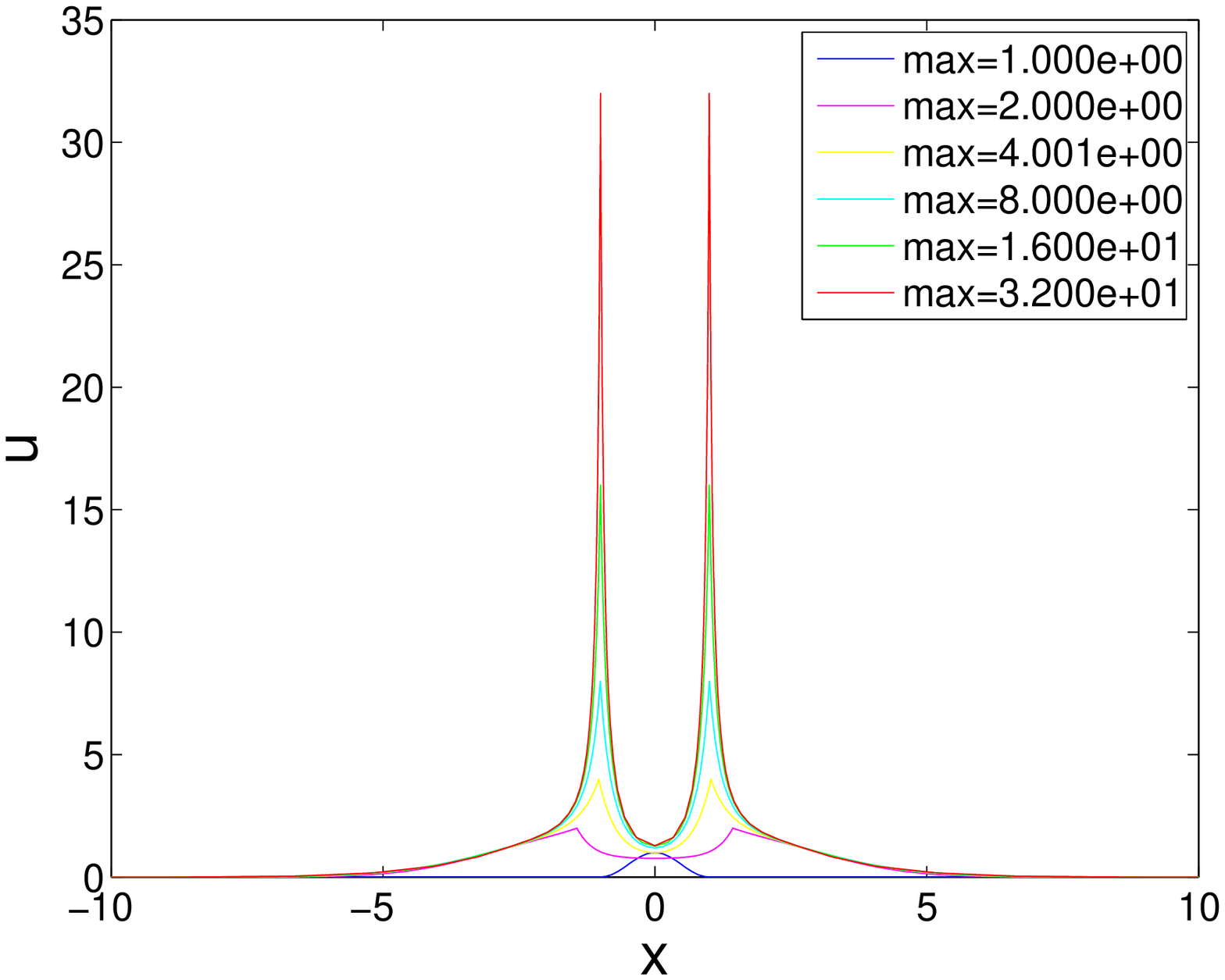}}%\\
  {\includegraphics[width=0.32\textwidth]{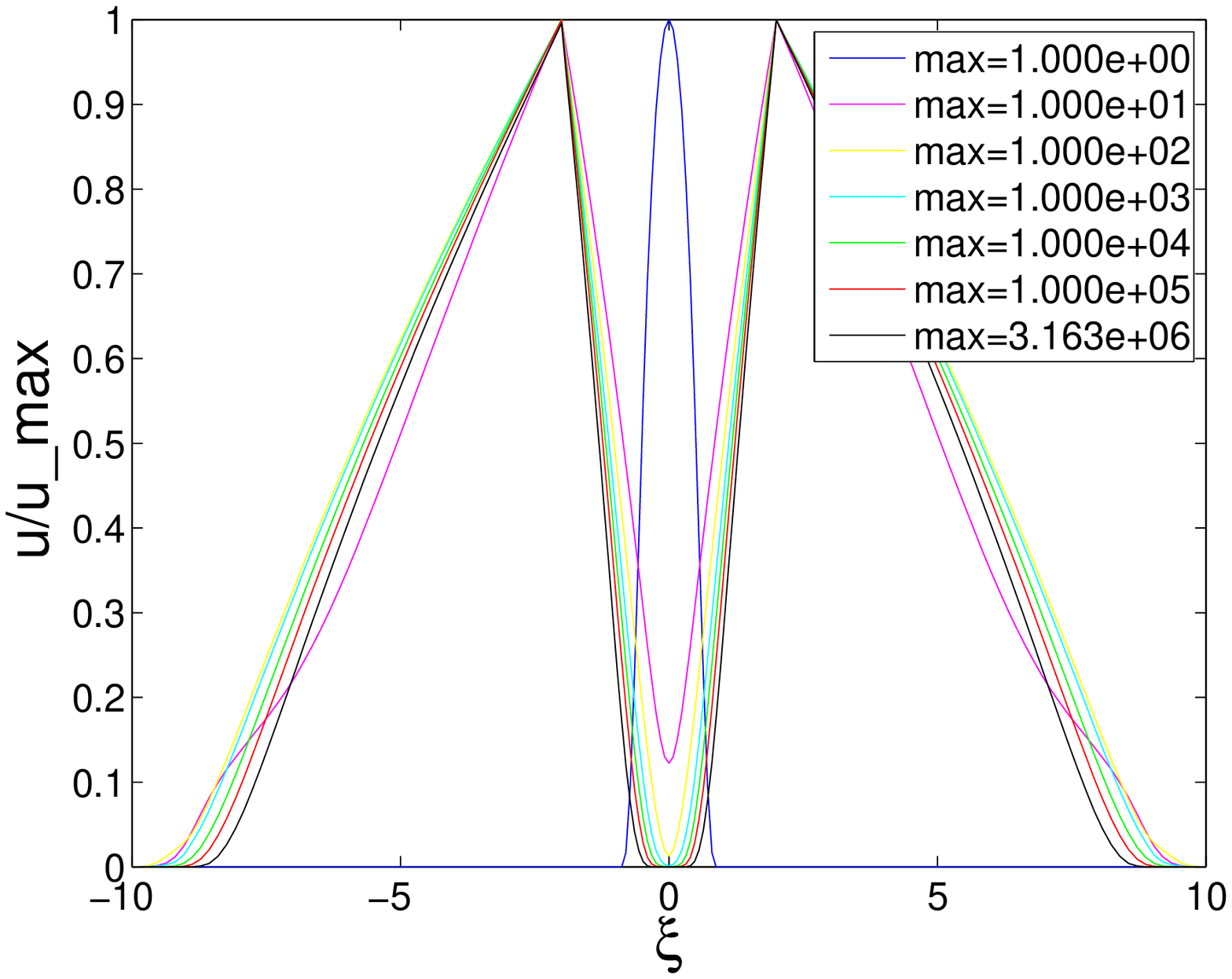}}%\hfill
  \caption{\small Numerical results for symmetric periodic moving sources with $N=150$.}
  \label{fig:sym1_source2}
\end{figure}

To our best knowledge, there has no theoretical results for multi-sources with different speeds and this is the first time numerically investigating the phenomenon for this case. It is shown in Figure \ref{fig:sym1_source2} that blow-up occurs on both sources at $t=2.496881990359248$, corresponding to the maximum value of $u_{\max}=3.16\times 10^6$.

\begin{figure}[!t]%[!htb]
  \centering
  {\includegraphics[width=0.36\textwidth]{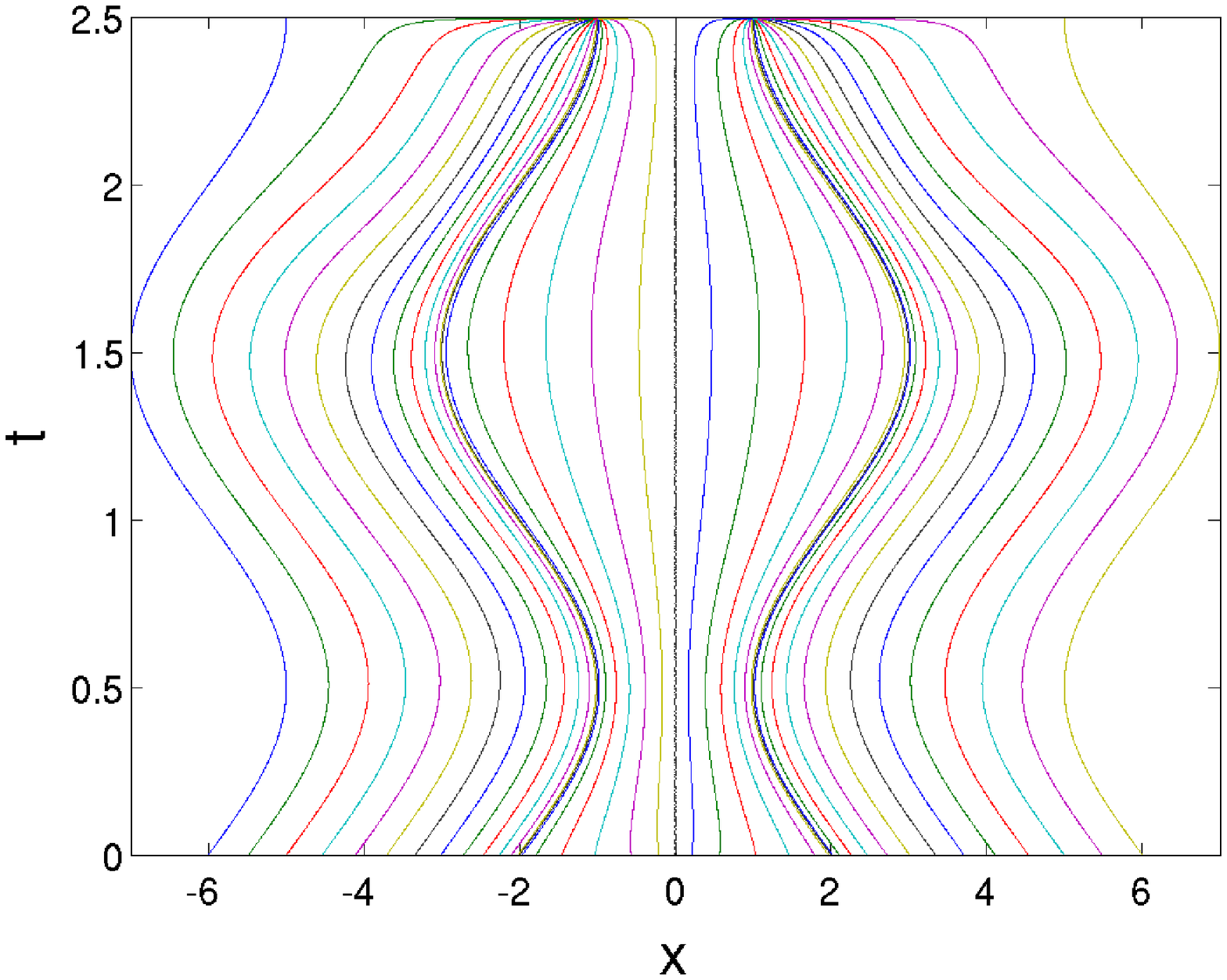}}%\hfill
  {\includegraphics[width=0.32\textwidth]{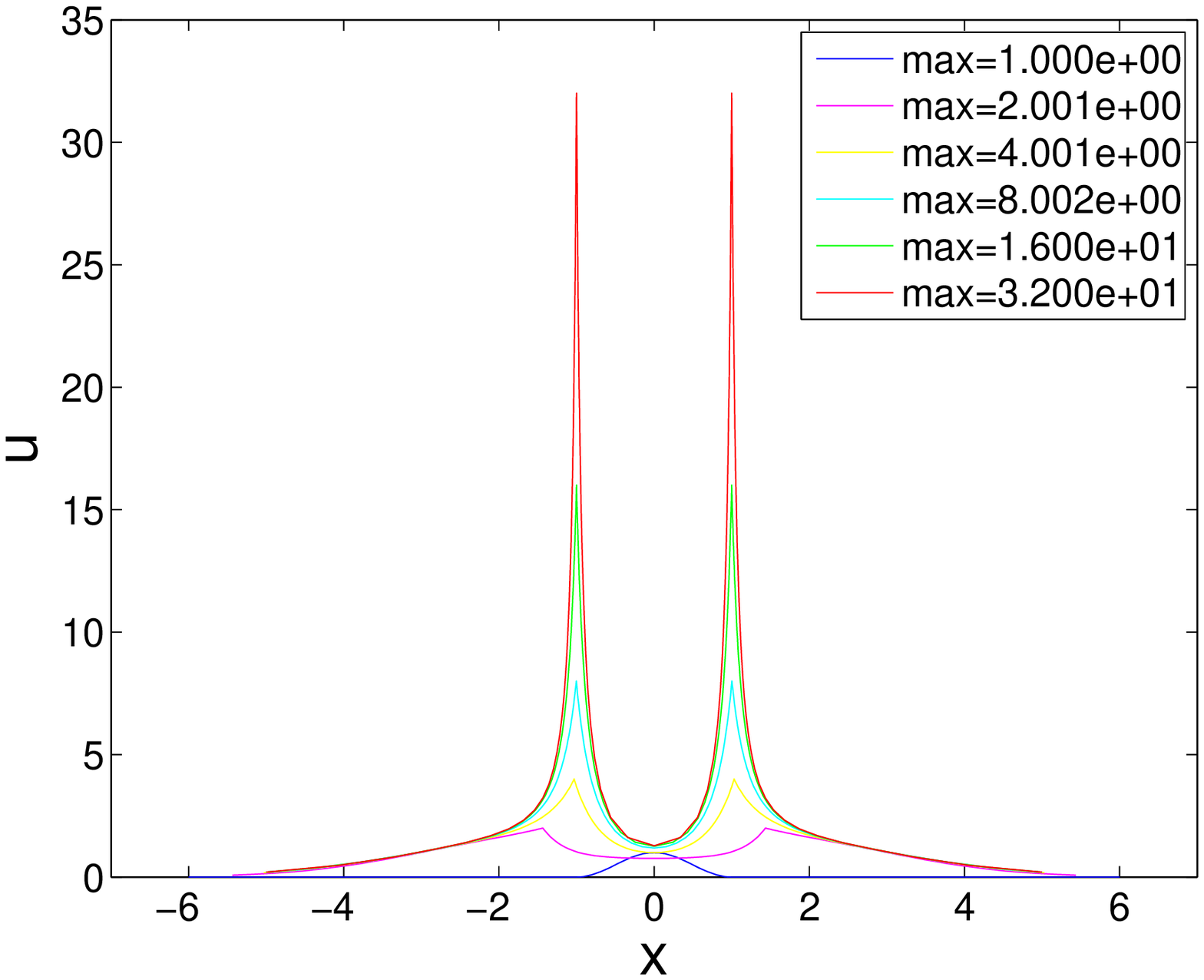}}%\\
  {\includegraphics[width=0.32\textwidth]{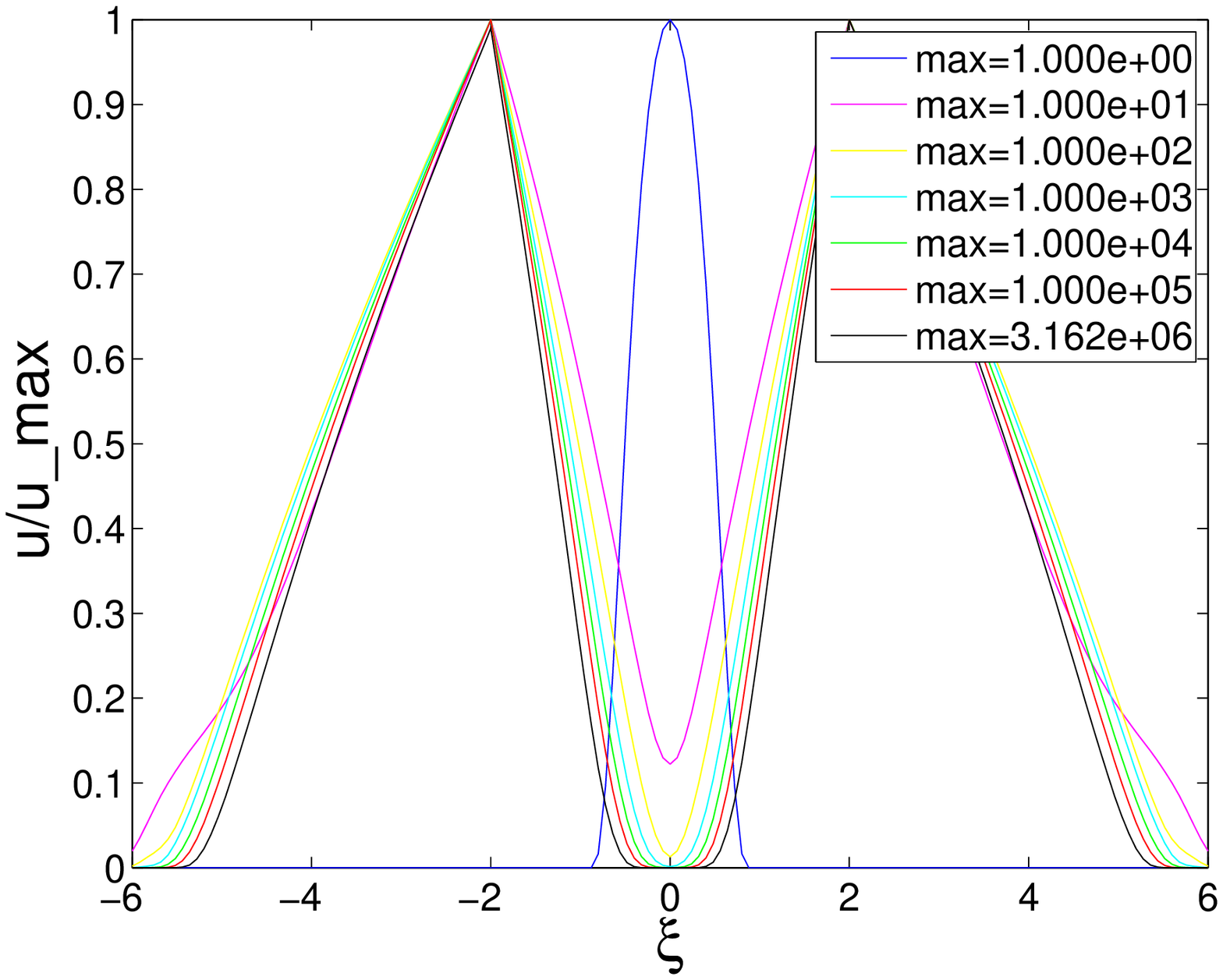}}%\hfill
  \caption{\small Numerical results for symmetric periodic moving sources with $[x_l,x_r]=[\alpha_0(t)-4.0,\alpha_1(t)+4.0]$.}
  \label{fig:sym1_source2_LABC}
\end{figure}

If local absorbing boundary conditions \eqref{eq:LABCs_log} are used, the observed domain can be chosen more smaller while the results do not be influenced almost. See Figure \ref{fig:sym1_source2_LABC} as an example, where the observed domain is set by $[\alpha_0(t)-4.0,\alpha_1(t)+4.0]$, changed as time evolution. Now blow-up occurs on both sources at $t=2.496370241342059$ with the maximum value of $u_{\max}=3.16\times 10^6$.

\section{Conclusions}\label{sec:conclusion}
In this paper, our work focus on the problem of traveling singular sources with different speeds.
A new moving mesh method in conjunction with a non-overlapping domain decomposition is proposed for solving this problems.
The whole domain is splitted into $q+1$ subdomains by the $q$ sources,
whose positions are gotten by a predictor-corrector algorithm.
Taking the advantages of the domain decomposition,
the computation of jump $[\dot{u}]$ is avoided and there are only two different cases discussed
in the discretization of the physical PDE.
Thus, it is easy for the implementation to solve the problems with two traveling sources or more.
Moreover, the moving mesh method of MMPDEs can be applied into each sub-domain respectively.
The second-order of the spatial convergence can be proved for the new method under
a special time marching implementation. The good performance of the new method for the blow-up phenomenon
is demonstrated through a number of examples with two sources.
%or more.
Furthermore, using the new method, we successfully simulate the solutions of two sources with different speeds. 
To our best knowledge, this is the first time investigation for this case. 
The case of three sources or more can be implemented similarly.

% In this paper, we have established a domain decomposed moving mesh method for traveling singular sources problems{\color{red}, where the sources' velocities are controlled by some ordinary differential equations 且不相交}. The domain is divided into $q+1$ subdomains by the $q$ sources, then the mesh on the observed domain is obtained by solving MMPDEs on each subdomain in parallel.
% % so that MMPDEs can be solved in parallel on each subdomain.
% Taking the advantages of that a fixed mesh point locates at each source during time evolution, the approximation scheme for the physical PDE is greatly simplified. A predictor-corrector algorithm is proposed to solve $u(x,t)$ and $\alpha_i(t)$, accurately. {\color{red}To our best knowledge, 还 没有任何 理论 结果, 在文献中仅仅考虑了同速运动的情形, 现有数值方法对于这种问题也不能有好的解决办法, 而我们能够有效解决这种问题的数值模拟, 我们的方法简单容易} Several numerical examples are presented to verify the convergence rates and efficiency of the method. Different type motion of the sources are considered for blow-up phenomenon. To our best knowledge, this is the first time investigating the phenomenon of the problem for multi-sources, which move with different speeds. Blow-up on both sources are also observed{\color{red}给理论分析提供了提示}. We will extend the method to interface problems and higher dimensional problems in the future{\color{red}相交情形}.

\section*{Acknowledgment}
This work was partially supported by a grant of key program from 
the National Natural Science Foundation of China (No.~10731060, 10801120, 11171305), 
National Basic Research Program of China (2011CB309704),  
Chinese Universities Scientific Fund No.~2010QNA3019 
and Zhejiang Provincial Natural Science Foundation of China under Grant No.~Y6110252.

% ===============================================
% =                参考文献                       =
% =           引文用 \cite{ 关键字 }               =
% ===============================================
\bibliographystyle{unsrt}
\bibliography{movingSource}

\end{document}